\documentclass[11pt]{amsart}
\usepackage{amsmath, amssymb}
 \usepackage{amsmath,amscd}
\usepackage{amsfonts}
\usepackage{mathrsfs}
\usepackage[arrow,matrix,curve,cmtip,ps]{xy}
\usepackage{enumitem}
\usepackage{setspace}

\usepackage{pb-diagram} 
\usepackage{pb-xy}
\usepackage{graphicx}

\usepackage{amsthm}

\usepackage{color}

\usepackage{xcolor}

\allowdisplaybreaks

\newtheorem*{rep@theorem}{\rep@title}
\newcommand{\newreptheorem}[2]{
\newenvironment{rep#1}[1]{
 \def\rep@title{#2 \ref{##1}}
 \begin{rep@theorem}}
 {\end{rep@theorem}}}
\makeatother

\newreptheorem{theorem}{Theorem}
\newreptheorem{lemma}{Lemma}
\newreptheorem{proposition}{Proposition}
\newreptheorem{corollary}{Corollary}

\newtheorem{theorem}{Theorem}[section]
\newtheorem{lemma}[theorem]{Lemma}

\newtheorem{proposition}[theorem]{Proposition}
\newtheorem{corollary}[theorem]{Corollary}

\newtheorem*{theorem*}{Theorem}
\newtheorem*{proposition*}{Proposition}
\theoremstyle{remark}
\newtheorem{remark}[theorem]{Remark}
\newtheorem{definition}[theorem]{Definition}

\numberwithin{equation}{section}

\newcommand{\Z}{\mathbb{Z}}
\newcommand{\Q}{\mathbb{Q}}

\newcommand{\R}{\mathbb{R}}
\newcommand{\Rn}{\mathcal{R}}
\newcommand{\C}{\mathcal{C}}

\newcommand{\K}{\mathcal{K}}
\newcommand{\F}{\mathcal{F}}
\newcommand{\cplx}{\mathbb{C}}

\newcommand{\im}{\operatorname{Im}}

\newcommand{\rank}{\operatorname{Rank}}

\newcommand{\widecap}{\widehat}
\newcommand{\infect}{{\operatorname{inf}}}
\newcommand{\into}{\hookrightarrow}
\newcommand{\onto}{\twoheadrightarrow}
\renewcommand{\hom}{\operatorname{Hom}}
\renewcommand{\ker}{\operatorname{Ker}}
\newcommand{\normalsubgroup}{\unlhd}

\newcommand{\nquotient}[2]{\ensuremath{\dfrac{\pi_1(#1)}{\pi_1(#1)_{\mathfrak{r}}^{(#2)}}}}

\newcommand{\iterate}[3]{\ensuremath{\underset{#2}{\overset{#3}{#1}}}}
\newcommand{\bdry}{\ensuremath{\partial}}

\newcommand{\eref}[1]{(\ref{#1})}

\begin{document}
\title[$\rho^1$ and infection.]{Linear Independence of Knots Arising from Iterated Infection Without the Use of Tristram Levine Signature.}

\author{Christopher Davis	}

\address{Department of Mathematics \\ Rice UNIVERSITY}
\email{cwd1@rice.edu}

\date{\today}



\begin{abstract}
We give an explicit construction of linearly independent families of knots arbitrarily deep in the ($n$)-solvable filtration of the knot concordance group using the $\rho^1$-invariant defined in \cite{myFirstPaper}.  A difference between previous constructions of infinite rank subgroups in the concordance group and ours is that the deepest infecting knots in the construction we present are allowed to have vanishing Tristram-Levine signatures.
\end{abstract}

\maketitle



\section{Introduction}
%

A knot $K$ is an isotopy class of oriented locally flat embeddings of the circle $S^1$ into the 3-sphere $S^3$.  A pair of knots $K$ and $J$ are called \textbf{topologically concordant} if there is a locally flat embedding of the annulus $S^1\times[0,1]$ into $S^3\times[0,1]$ mapping $S^1\times\{1\}$ to a representative of $K$ in $S^3\times\{1\}$ and $S^1\times\{0\}$ to a representative of $J$ in $S^3\times\{0\}$.  A knot is called \textbf{slice} if it is concordant to the unknot or equivalently if it is the boundary of a locally flat embedding of the 2-ball $B^2$ into the 4-ball $B^4$.  The set of all knots modulo concordance under the operation of connected sum is a group called the \textbf{knot concordance group} and is denoted by $\C$.

In \cite{whitneytowers}, Cochran, Orr and Teichner define the solvable filtration of $\C$:
\begin{equation*}
\dots \F_{n.5} \subseteq \F_{n} \subseteq \dots \subseteq \F_{1.5}\subseteq \F_1\subseteq \F_{0.5}\subseteq \F_0\subseteq \C.
\end{equation*}
For $k$ a half integer, the elements in $\F_k$ are called \textbf{($k$)-solvable}.  They show that $\F_0$ is the set of Arf-invariant zero knots, $\F_{0.5}$ is the set of algebraically slice knots and that Casson-Gordon invariants vanish on $\F_{1.5}$.  In \cite[Section 6]{whitneytowers} Cochran-Orr-Teichner show that $\F_2/\F_{2.5}$ is infinite rank by studying a satellite operation (called infection in \cite[Section 8]{C}).  The quotient groups $\F_n/\F_{n.5}$ have been an active place of research ever since.  In \cite{primaryDecomposition} Cochran, Harvey and Leidy begin with knots for which the integrals of the Tristram-Levine signature functions are linearly independent over $\Q$ and use an iterated infection procedure to produce an infinite rank free Abelian subgroup of $\F_n/\F_{n.5}$.

In \cite{amenableL2Methods} Cha constructs another infinite rank subgroup of $\F_n/\F_{n.5}$ starting instead with knots whose Tristram-Levine signatures evaluate to sufficiently large values at particular finite sets.  

We present a variation on this idea, performing iterated infections to produce linearly independent sets deep in the solvable filtration.  A novel aspect of the construction we present is that the deepest infecting knots are allowed have vanishing Tristram-Levine signature.  A concrete advantage of this construction over previous ones is that its conditions may be directly verified, providing explicit infinite linearly independent sets rather than generating sets for infinite rank subgroups.

Given any closed oriented 3-manifold $M$ and a homomorphism $\phi:\pi_1(M) \to \Gamma$, the von Neumann $\rho$-invariant, $\rho(M, \phi)\in \R$, is defined.  It is an invariant of orientation preserving homeomorphism of the pair $(M,\phi)$.  Restricting this invariant to the zero surgery of knots and links gives rise to an isotopy invariant.  We provide a brief overview of $\rho$-invariants in Section~\ref{background}.

In \cite{myFirstPaper} the author defines a particular $\rho$-invariant, $\rho^1$, shows that in a restricted setting it provides a concordance obstruction and computes it for an infinite family of twist knots of order 2 in the algebraic concordance group.  In this paper we bring this invariant to bear on an iterated infection procedure in order to produce examples whose deepest infecting knots have vanishing Tristram-Levine signature.

For an overview of infection, see \cite[section 8]{C}.  We denote the infection of the base knot $R$ along the infecting curve $\eta$ in $S^3-R$ by the infecting knot $J$ as either $R_\eta(J)$ or $R(\eta,J)$ depending on notational convenience.

Recall that for a knot $K$, the \textbf{rational Alexander module} of $K$, $A_0(K)$, is given by the first homology with coefficients in $\Q$ of the infinite cyclic cover of the exterior of $K$ or equivalently of the zero surgery of $K$, $M(K)$.  The rational Alexander module of $K$ is a module over the ring of Laurent polynomials, $\Q[t^{\pm1}]$.  With respect to the involution on $\Q[t^{\pm1}]$ given by $\overline{p(t)}=p(t^{-1})$, there is a sesquilinear form $$Bl:A_0(K)\times A_0(K)\to \frac{\Q(t)}{\Q[t^{\pm1}]},$$ called the \textbf{Blanchfield form}.  A submodule $P\subseteq A_0(K)$ is called \textbf{isotropic} if $Bl(x,y)$ vanishes for all $x,y\in P$.

We now give definitions of the concepts needed in the statement of Theorem~\ref{big corollary}, the main theorem of this paper.  They will be recalled when needed.  

Two polynomials $p(t), q(t)$ are called \textbf{strongly coprime} if $p(t^k)$ and $q(t^l)$ have no common roots in $\cplx$ for every choice of nonzero integers $k$ and $l$.  A pair $(R,\eta)$ with $R$ a knot and $\eta$ a curve in its complement is called \textbf{doubly anisotropic} if $\eta$ represents an element of  $A_0(R)$  for which there does not exist any $\alpha$ and $\beta$ in $A_0(R)$ with $\eta=\alpha+\beta$ and $Bl(\alpha,\alpha)=Bl(\beta,\beta)=0$.  

\begin{theorem}\label{big corollary}
Let $\{K_i\}$ be a possibly infinite set of knots:
\begin{enumerate}
\item whose Alexander polynomials are strongly coprime,
\item whose Tristram-Levine signatures have vanishing integrals,
\item whose prime factors have square-free Alexander polynomials and
\item whose $\rho^1$-invariants do not vanish, that is $\rho^1(K_i)\neq 0$.
\end{enumerate}
   For $i=1,2, \dots$ and $j=1,2,\dots, n$ let $R_{i,j}$ be a slice knot and $\eta_{i,j}$ be an unknotted curve in the complement of $R_{i,j}$ such that the pair $(R_{i,j},\eta_{i,j})$ is doubly anisotropic.  Let $K_i^0=K_i$ and $K_i^j=R_{i,j}(\eta_{i,j},K_i^{j-1})$.

 Then $\left\{K_i^n\right\}_{i=1}^\infty$ is linearly independent in $\C$ modulo ($n+1.5$) solvable knots.
\end{theorem}

A noteworthy difference between this and previous results producing infinite rank subgroups via iterated infection is the condition on the Tristram-Levine signature.  Previous constructions assume that the deepest infecting knots, $K_{i}$, be complicated in some sense.  In \cite[Theorem 7.5]{primaryDecomposition} their integrals are required to be rationally linearly independent.  In \cite[Lemma 4.11 and Proposition 4.12]{amenableL2Methods} they are required to take large values at a specific set of points.  By contrast, Theorem~\ref{big corollary} is designed to apply even when Levine-Tristram signature functions vanish. 

Another key difference between this result and such previous results is the ease of verifying that the assumptions of the theorem are satisfied.  Specifically, notice that in Theorem~\ref{big corollary} there is no assumption on the $\rho$-invariants of the slice knots $R_{i,j}$.  The techniques of  \cite{primaryDecomposition} require a condition on first order von Neumann $\rho$-invariants of the knots along which infection is performed.  Without a means of computation, they cannot verify that any fixed knot satisfies this condition.  The techniques of \cite{amenableL2Methods} require that the Tristram-Levine signature of the infecting knots exceed all of the von Neumann $\rho$-invariants of the knots along which infection is done.  Without any means of getting concrete bounds, these techniques will not give any explicit linearly independent sets.

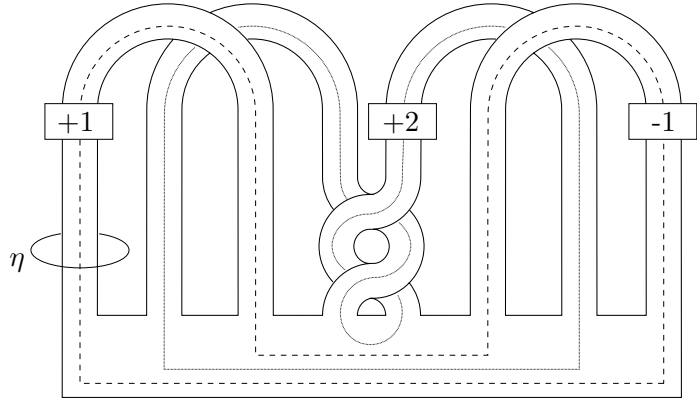
\begin{figure}[htbp]
\setlength{\unitlength}{1pt}
\begin{picture}(150,150)
\put(-45,0){\includegraphics[width=.7\textwidth]{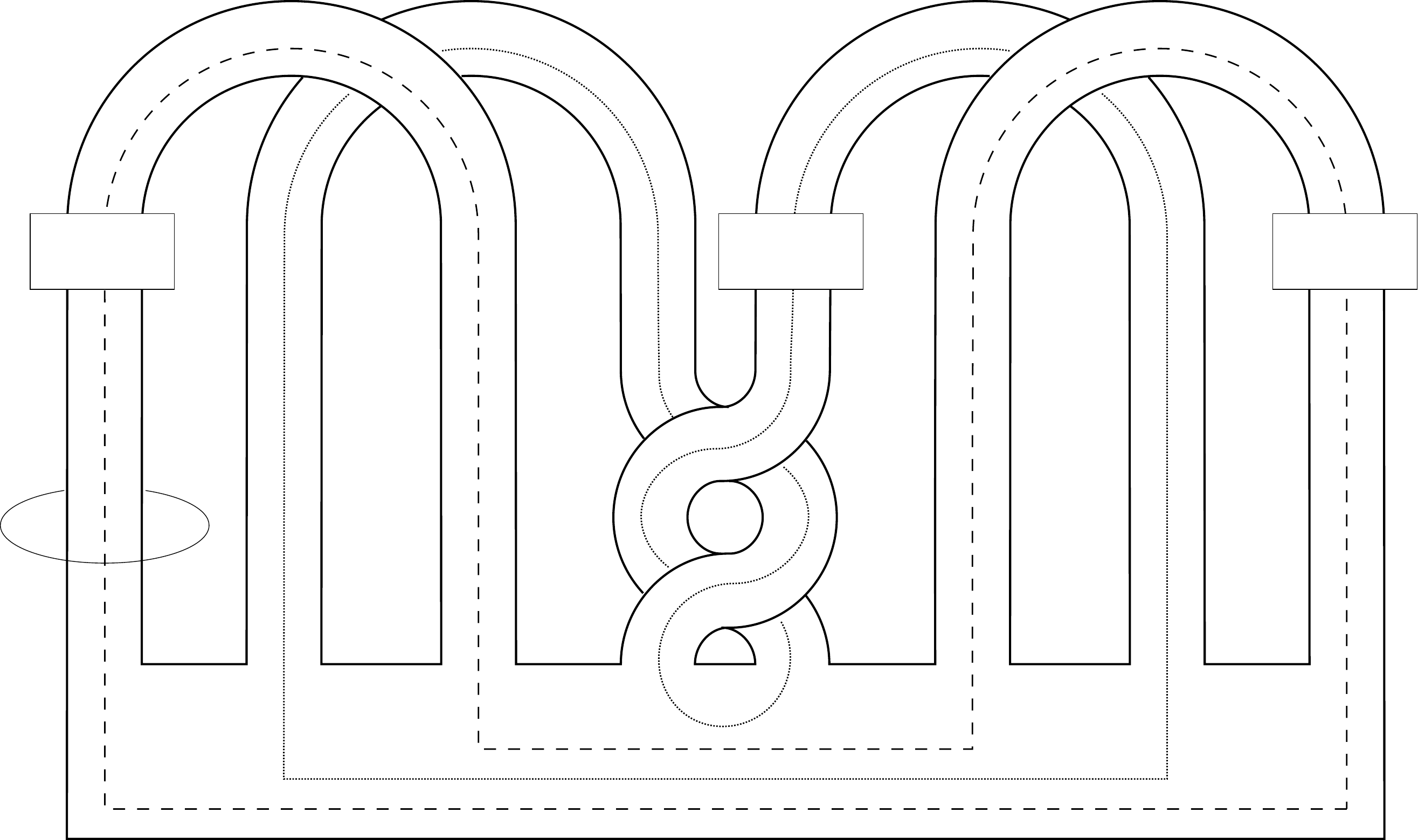}}
\put(-35,101){+1}
\put(190,101){-1}
\put(88,101){+2}
\put(-53,50){$\eta$}
\end{picture}
\caption{A slice knot $R$ with a doubly anisotropic curve, $\eta$.  The depicted derivative is the unlink.}\label{fig:slice example}
\end{figure}

\begin{figure}[htbp]
\setlength{\unitlength}{1pt}
\begin{picture}(80,100)
\put(-20,0){\includegraphics[width=.3\textwidth]{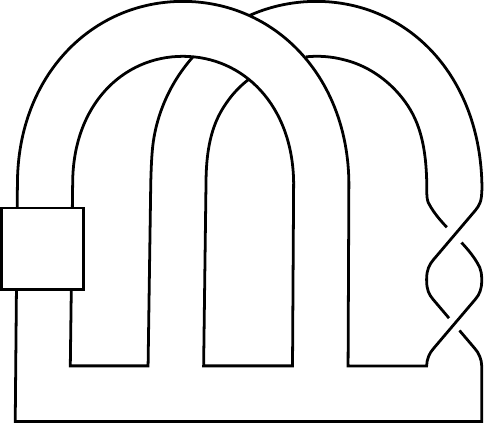}}
\put(-15,35){\Large{$n$}}
\end{picture}
\caption{$T_n$, the $n$-twist knot.}\label{fig:twistknot}
\end{figure}

 In Section~\ref{examples}, as an application of Theorem~\ref{big corollary} we generate the following family of linearly independent knots deep in the filtration using the base knot $R$, infecting curve $\eta$ (see Figure~\ref{fig:slice example}) and deepest infecting knot given by twist knots $T_n$ of finite algebraic order (see Figure~\ref{fig:twistknot}).  

\begin{reptheorem}{example}
For the slice knot $R$ and infecting curve $\eta$, if $T_n$ is the $n$-twist knot then
\begin{equation*}
\{(R_\eta)^m(T_n) = (R_\eta\circ\dots\circ R_\eta)(T_n)|n=-x^2-x-1, x\ge2\}
\end{equation*}
is linearly independent in $\F_{m-.5}/\F_{m+1.5}$, where $\F_{-.5}$ is taken to be the whole concordance group.
\end{reptheorem}

This family of knots appears to be the first linearly independent set deep in the solvable filtration of $\C$ constructed by an iterated infection procedure with deepest infecting knots whose Tristram-Levine signature functions vanish.

\subsection{Outline of the paper}

In Section~\ref{background} we provide what in this paper is taken as the definition of the von Neumann $\rho$-invariant, as well as some properties of the $L^2$-signature.  We go on to provide definitions of the $\rho^0$ and $\rho^1$-invariants which are employed in this paper.  

In Section~\ref{examples}, we find a family of knots with nonzero $\rho^1$-invariant whose Alexander polynomials are strongly coprime and square-free, as well as a set of slice knots whose rational Alexander modules have doubly anisotropic elements, that is, a set of knots satisfying the conditions of  Theorem~\ref{big corollary}.  This provides an explicit  linearly independent set of knots sitting arbitrarily deep in the solvable filtration of the concordance group.  The remaining sections are devoted to the development of the machinery used in the proof of Theorem~\ref{big corollary}

In Section~\ref{invariant} we discuss a localization of the Alexander module, $\widetilde{A_0^p}(K)$.  In order to capture information involving this localization we define a new class of von Neumann $\rho$-invariant.  It enjoys additivity properties over connected sum and infection and in some cases agrees with $\rho^1$.

In Section~\ref{anisotropy examples} we study the Blanchfield linking form on $\widetilde{A_0^p}(K)$ and find sufficient conditions for this localized Blanchfield form to have no nontrivial isotropic submodules.  The significance of this result appears in Section~\ref{isotropy} in which we show that isotropic submodules of this localization are well behaved with respect to the operation of extension of scalars.   

 Finally, in Section~\ref{main section} we give the proof of Theorem~\ref{big corollary}.

\section{Background: von Neumann $\rho$-invariants and $L^2$ signatures}\label{background}

In this section we state the properties of von Neumann $\rho$-invariants and $L^2$ signatures needed in this paper.  

In \cite[Section 3]{Ha2} the following property is proven of the von Neumann $\rho$-invariant.  It serves here as the definition.  

\begin{definition}\label{rho}
Consider oriented 3-manifolds $M_1, \dots, M_n$, with homomorphisms $\phi_i:\pi_1(M_i)\to \Gamma_i$.  Suppose that $M_1\sqcup M_2\sqcup\dots\sqcup M_n$ is the oriented boundary of a compact oriented 4-manifold $W$ and $\psi:\pi_1(W)\to \Lambda$ is a homomorphism such that, for each $i$, there is a monomorphism $\alpha_i:\Gamma_i\to \Lambda$ making the following diagram commute:
\begin{center}$
 \begin{diagram}
\node{\pi_1(M_i)} \arrow{e,t}{\phi_i}
         \arrow{s,r}{i_*}
	     \node{\Gamma_i} \arrow{s,r,J}{\alpha_i}\\
       \node{\pi_1(W)} \arrow{e,t}{\psi}
     \node{\Lambda} 
      \end{diagram}$
      \end{center}

Then $\displaystyle\sum_{i=1}^n\rho(M_i,\phi_i) = \sigma^{(2)}(W,\psi) - \sigma(W)$ where $\sigma(W)$ is the regular signature of $W$ and $\sigma^{(2)}(W,\psi)$ is the $L^2$ signature of $W$ twisted by the coefficient system $\psi$.  The expression $\sigma^{(2)}(W,\psi) - \sigma(W)$ is called the signature defect of $W$ with respect to $\psi$.
\end{definition}

For a compact oriented 4-manifold $W$ with coefficient system $\phi:\pi_1(W)\to \Gamma$, $\sigma^{(2)}(W,\Gamma)\in \R$ is defined.  The properties of the $L^2$ signature which are used in this paper are Novikov additivity and a bound in terms of ranks of twisted second homology in the case that $\Gamma$ is PTFA (Poly Torsion Free Abelian, see \cite[Definition 2.1]{whitneytowers}).

\begin{proposition*}[{Novikov additivity, \cite[Lemma 5.9 (3)]{whitneytowers}}] \label{novikov}
If compact oriented 4-manifolds $W_1$ and $W_2$ intersect in a single common boundary component, $W=W_1\cup W_2$, and $i^1:W_1\to W$ and $i^2:W_2\to W$ are the inclusion maps, then for every homomorphism $\phi:\pi_1(W)\to \Gamma$, $\sigma^{(2)}(W,\phi) = \sigma^{(2)}(W_1,\phi\circ i^1_*)+\sigma^{(2)}(W_2,\phi\circ i^2_*)$.
\end{proposition*}

The second property is that when $\Gamma$ is PTFA and more generally whenever $\Q[\Gamma]$ is an Ore domain,
\begin{equation}\label{inequality}
\left|\sigma^{(2)}(W,\phi)\right| \le \rank_{\Q[\Gamma]}\left(\dfrac{H_2\left(W;\Q[\Gamma]\right)}{i_*\left[H_2\left(\bdry W;\Q[\Gamma]\right)\right]}\right)
\end{equation}
where $i_*:H_2\left(\bdry W;\Q[\Gamma]\right)\to H_2\left(W;\Q[\Gamma]\right)$ is the inclusion induced map.  This follows from the monotonicity of von Neumann dimension (see \cite[Lemma 1.4]{L2invts}) and the fact that the $L^2$ Betti number agrees with $\Q[\Gamma]$ rank when $\Q[\Gamma]$ is an Ore Domain (see \cite [Lemma 2.4] {Cha3} or \cite[Proposition 2.4]{FrLM}). 

\subsection{The $\rho^0$ and $\rho^1$-invariants}
\begin{definition}
For a knot $K$, Let $\phi^0:\pi_1(M(K))\to \Z$ be the Abelianization map.  Let $\rho^0(K) := \rho(M(K),\phi^0)$ be the corresponding $\rho$-invariant.
\end{definition}

It is shown in \cite[Proposition 5.1]{structureInConcordance} that $\rho^0(K)$ is given by the integral of the Tristram-Levine signature function.  

Recall that the rational derived series of a group $G$ is given by setting  $G^{(0)}_\mathfrak{r}=G$ and inductively defining $G^{(n+1)}_\mathfrak{r}$ to be the set of all $g\in G^{(n)}_\mathfrak{r}$ which are torsion in the Abelianization of $G^{(n)}_\mathfrak{r}$.  This series is the most quickly descending series with the property that each of the successive quotients $G^{(n)}_\mathfrak{r}/G^{(n+1)}_\mathfrak{r}$ is TFA (Torsion Free Abelian).

\begin{definition}
For a knot $K$, let \begin{equation*}\phi^1:\pi_1(M(K))\to \nquotient{M(K)}{2}\end{equation*} be the quotient by the second term in the rational derived series.  Let $\rho^1(K) := \rho(M(K),\phi^1)$ be the corresponding $\rho$-invariant.
\end{definition}

In \cite{myFirstPaper} the $\rho^1$-invariant is shown to provide a sliceness obstruction and is used to find an infinite collection of twist knots of algebraic order 2 which is linearly independent in $\C$.

\section{generating explicit linearly independent families of knots arbitrarily deep in the solvable filtration}\label{examples}

In this section we verify the assumptions of Theorem~\ref{big corollary} for an explicit set of knots.     

We first address the restriction we put on the slice knots we infect and curves along which we infect them.

\begin{definition}

For a knot, $R$, an element of $A_0(R)$, $\eta$, is called \textbf{doubly anisotropic} if $\eta$ cannot be written as a sum of isotropic elements, that is there do not exist any $\alpha, \beta\in A_0(R)$ with $Bl(\alpha,\alpha)=Bl(\beta,\beta)=0$ and $\alpha+\beta=\eta$.  
\end{definition}

In this paper we concern ourselves with the operator $R_\eta:\C\to\C$ given by sending a knot $J$ to the result of infection $R_\eta(J)$ where $R$ is a slice knot and $\eta$ is an unknotted curve representing a doubly anisotropic element of $A_0(R)$. 

The following proposition serves to illustrate that there are many slice knots whose Alexander modules have doubly anisotropic elements.  

\begin{proposition}\label{doubly anisotropic example}
Let $K$ be a slice knot with cyclic Alexander module isomorphic to $\frac{\Q[t^{\pm1}]}{\left(\delta(t)^2\right)}$ where $\delta$ is any symmetric polynomial with $\delta(1)=\pm1$.  If $\delta$ has a prime symmetric factor then $\eta$, the generator of $A_0(R)$, is doubly anisotropic.  
\end{proposition}

\begin{proof}

Let $q$ be the assumed prime symmetric factor of $\delta$.  Consider any element, $f\eta$, in the Alexander module ($f\in \Q[t^{\pm1}]$).  If $f\eta$ were an isotropic element then \begin{equation*}0=Bl(f\eta,f\eta)=\displaystyle \dfrac{f(t)f(t^{-1})r(t)}{\delta(t)^2}\in \dfrac{\Q(t)}{\Q[t^{\pm1}]}\end{equation*}
 where $(r,\delta)=1$.  Then $q$, being a factor of $\delta$, must divide $f(t)f(t^{-1})r(t)$.  Since $q$ is prime and $(r,q)=1$ $q$ must divide $f(t)$ or $f(t^{-1})$.  Since $q$ is symmetric, it divides both.  Thus, any isotropic element of the $A_0(R)$ (and so any sum of two isotropic elements) sits in the proper submodule $P=\langle q\eta\rangle$ and any element of $A_0(R)-P$ (for example, $\eta$, the generator of $A_0(R)$) is doubly anisotropic.
\end{proof}

Cha (\cite[Theorem 5.18]{Cha2}) shows that for every symmetric polynomial $\delta$ there is a ribbon knot with Alexander module of the form $\frac{\Q[t^{\pm1}]}{(\delta(t)^2)}$, so that slice knots with doubly anisotropic curves abound.  For the sake of concreteness, let $R$ be the slice knot depicted in Figure~\ref{fig:slice example} and Let $\eta$ be the curve in $S^3-R$, also depicted in Figure~\ref{fig:slice example}.  The Alexander module of $R$ is cyclic generated by $\eta$ and has Alexander polynomial of the form $\delta(t)^2$ where $\delta(t)=t^2-3t+1$ is a symmetric prime polynomial.  By Proposition~\ref{doubly anisotropic example} $(R,\eta)$ is doubly anisotropic.

What remains is to find infinitely many infecting knots whose Tristram-Levine signatures have vanishing integrals, whose $\rho^1$-invariants are nonzero and whose Alexander polynomials are strongly coprime and square-free.  For $n<0$, the twist knot $T_n$ (see Figure~\ref{fig:twistknot}) is algebraically of finite order, so that the Tristram-Levine signature vanishes.  It is shown in \cite[Theorem 6.1]{myFirstPaper} that for $n(x)=-x^2-x-1$, and $x\ge2$ $\rho^1(T_{n(x)})\neq 0$.  Their Alexander polynomials are prime and so square-free.  

Theorem 3.1 of \cite{paperWithEvan} gives us that the strong coprimality condition is satisfied.  A note on conventions, the knot which is called $T_n$ in this paper is the reverse of the mirror image of the knot called $T_{-n}$ in \cite{paperWithEvan}.

\begin{theorem*}[Theorem 3.1, \cite{paperWithEvan}]
For positive integers $m\neq n$  the Alexander polynomials $\Delta_{T_{-n}}$ and $\Delta_{T_{-m}}$ are strongly coprime.
\end{theorem*}

Thus, the slice knot $R$ and infecting curve $\eta$, together with the deepest infecting knot $T_{n(x)}$ for $x\ge 2$ satisfy the assumptions of Theorem~\ref{big corollary} and we see that

\begin{theorem}\label{example}
The set $\{(R_\eta)^m(T_{n})|n=-x^2-x-1,x\ge 2\}$ is linearly independent in $\F_{m-0.5}/\F_{m+1.5}$, where $\F_{-0.5}$ is taken to be all of the concordance group.
\end{theorem}

Starting with a different choice of $R_\eta$, we construct families of knots that are linearly independent in the concordance group but which many previous invariants fail to detect.

\begin{theorem}
Let $p_m(t)$ denote the $m$th cyclotomic polynomial where $m$ is divisible by three distinct prime  numbers.

Let $R$ be a ribbon knot with cyclic Alexander module $A_0(R)\cong\frac{\Q[t^{\pm1}]}{(p_m^2)}$.  Let $\eta$ be an unknotted curve representing a generator of $A_0(R)$.     Let $T_n$ be the $n$-twist knot.  

Then $\{R_\eta(T_n)|n=-x^2-x-1\text{ with } x\ge 2\}$ is linearly independent in $\F_{0.5}/\F_{2.5}$ (and so in $\C$); however, the Casson-Gordon sliceness obstruction of \cite{CG2}, the metabelian $\eta$-invariant obstruction of \cite{fr3} and the $(1.5)$-solvability obstructions of \cite{polynomialSplittingOfRho}, \cite{whitneytowers} and \cite{amenableL2Methods}  all vanish for each element of this set.
\end{theorem}
\begin{proof}

The fact that this set is linearly independent is an immediate consequence of Theorem~\ref{big corollary}.

For every $n$ the Alexander polynomial of $R_\eta(T_n)$ is the same as the Alexander polynomial of $R$, which is $p_m^2$.  By \cite[Theorem 1.2]{Li7}, every prime power cyclic branched cover of $R_\eta(T_n)$ is a homology sphere.  Thus, the metabelian $\eta$-invariants of \cite{fr3} and the Casson-Gordon obstructions vanish.

In order to compute the obstructions of \cite{whitneytowers} (using the specialization of the obstruction to ($1.5$)-solvability of \cite[Theorem 4.2]{derivatives}) and \cite[Theorem 4.1]{polynomialSplittingOfRho}, notice that there is only one Lagrangian submodule of $A_0(R_\eta(T_n))\cong A_0(R)$, namely $P=\langle p_m\rangle$.  We first compute the obstruction of \cite{derivatives}.   By \cite[Lemma 2.3]{blanchfieldDuality}
\begin{eqnarray*}
\rho^1_{P}(R_\eta(T_n))&=&\rho^1_{P}(R)+\rho^0(T_n).
\end{eqnarray*}
As $T_n$ is of finite algebraic order, $\rho^0(T_n)=0$.  Since $R$ is slice and $P$ is the only Lagrangian submodule of $A_0(R)$, \cite[Theorem 4.2]{derivatives} implies that $\rho^1_{P}(R)=0$.  This completes the proof that this invariant vanishes.

In order to compute the obstruction of \cite{polynomialSplittingOfRho}, let $x\in P$ and consider the map $\phi_x$ defined in that paper.  By \cite[Lemma 2.3]{blanchfieldDuality}, $$\rho(R_\eta(T_n), \phi_x)=\rho(R,\phi_x)+\rho^0(T_n),$$ similarly to before, $\rho^0(T_n)$ vanishes and since $R$ is slice and has only one Lagrangian submodule, $P$, \cite[Theorem 1.1]{polynomialSplittingOfRho} implies that $\rho(R,\phi_x)=0$ for all $x\in P$.

Finally, we check that the $(1.5)$-solvability obstruction of \cite[Theorem 1.3]{amenableL2Methods} vanishes.  Since $R$ is slice, it follows that  $\rho(M(R),\phi)=0$ for some coefficient system $\phi:\pi_1(M(R))\to \Gamma$ where $\Gamma^{(2)}=0$, $\Gamma$ is amenable and $\Gamma$ is in Strebel's class $D(S)$ (See \cite{amenableL2Methods} for a definition) where $S$ is either $\Q$ or a finite cyclic group.  By \cite[Lemma 2.3]{blanchfieldDuality}, then $$\rho(M(R_\eta(T_n)),\phi)=\rho(M(R),\phi)+\rho(M(T_n),\phi).$$
As we observed, $\rho(M(R),\phi)=0$.  By \cite[Lemma 4.5]{amenableL2Methods}, $\rho(M(T_n),\phi)$ is a sum or integral of the Tristram-Levine signature of $T_n$ and so is zero.

\end{proof}

We do not know if it is possible to use the ($2.5$)-solvability obstructions of \cite{whitneytowers} or \cite{amenableL2Methods} to show that these knots are not ($2.5$)-solvable.

\section{strongly localized $\rho$-invariants}\label{invariant}

For a knot $K$ and a polynomial $p$, the strongly localized $\rho$-invariant of $K$, $\widetilde{\rho^1_p}(K)$ is defined in terms of a localization of the Alexander module of $K$.  We begin by describing this localization.  

For polynomials $p,q\in \Q[t^{\pm1}]$ we say that $p$ and $q$ are \textbf{strongly coprime} (\cite[Definition 4.4]{primaryDecomposition}) denoted $\widetilde{(p,q)}=1$ if there is no nonzero complex number $z$ and integers $m,n$ such that $p(z^m)=q(z^n)=0$.  Let \begin{equation*}\widetilde{S_p}=\{q\in\Q[t^{\pm1}] | \widetilde{(p,q)}=1\} \end{equation*} be the multiplicative set consisting of polynomials strongly coprime to $p$.

Let 
\begin{equation*}
{\widetilde{R_p}} 
=
 \Q[t^{\pm1}]\widetilde{S_p}^{-1}
=
\left\{\left.\dfrac{f}{g}\in\Q(t)\right|\widetilde{(g,p)}=1\right\}
\end{equation*}
 be the \textbf{strong localization} of $\Q[t^{\pm1}]$ at $p$.  By \cite[Theorem 10.30]{Rotman} $\widetilde{R_p}$ is flat as a $\Q[t^{\pm1}]$ module so that for a knot $K$, the first homology of $M(K)$ with coefficients in $\widetilde{R_{p}}$ is given by $H_1(M(K);\widetilde{R_p})\cong H_1(M(K);\Q[t^{\pm1}])\underset{\Q[t^{\pm1}]}{\otimes}\widetilde{R_p}$.  We call this module the \textbf{strongly localized Alexander module} and denote it by $\widetilde{A_0^p}(K)$.

Let $\pi_1(M(K))^{(2)}_{\widetilde p}$ be the kernel of the composition
\begin{equation}\label{composition}\pi_1(M(K))^{(1)}\to \frac{\pi_1(M(K))^{(1)}}{\pi_1(M(K))^{(2)}}\into A_0(K)\to \widetilde{A_0^p}(K).\end{equation}
Let 
$\widetilde{\phi_1^p}:\pi_1(M(K))\to \frac{\pi_1(M(K))}{\pi_1(M(K))^{(2)}_{\tilde{p}}}$ be the quotient map.  
\begin{definition}
Let ${\widetilde{\rho^1_p}}(K)=\rho(M(K),\widetilde{\phi_1^p})$ be the \textbf{strongly localized $\rho$-invariant of $K$ at $p$}.
\end{definition}

We will be flexible with notation.  For any CW-complex $X$ with infinite cyclic first homology generated by $t$ we can similarly define ${A_0}(X)$, $\widetilde{A_0^p}(X)$ and $\pi_1(X)^{(2)}_{\widetilde p}$.

The strongly localized $\rho$-invariant shares many properties with the similarly defined localized $\rho$-invariant of \cite{myFirstPaper}.  The proofs of the following propositions are identical to proofs in \cite{myFirstPaper} and so are omitted.  

\begin{proposition}[{\cite[Proposition 3.4]{myFirstPaper}}]\label{easy}  

If $\Delta$ is the Alexander polynomial of a knot $K$, then 
\begin{enumerate}
\item
$\widetilde{\rho^1_\Delta}(K)=\rho^1(K)$
\item
If $p$ and $\Delta$ are strongly coprime, then $\widetilde{\rho^1_p}(K)=\rho^0(K)$.
\end{enumerate}
\end{proposition}
%
%
%

\begin{proposition}[{\cite[Proposition 3.5]{myFirstPaper}}]\label{homomorphism}

  Let $J$ and $K$ be knots and $\eta$ be an unknot in the complement of $J$ such than $J$ and $\eta$ have zero linking number.

\begin{enumerate}
\item $\widetilde{\rho^1_p}(J\# K) = \widetilde{\rho^1_p}(J) + \widetilde{\rho^1_p}(K)$
\item $\displaystyle \widetilde{\rho^1_p}(J_\eta(K)) = \left\{ 
\begin{array}{ccc}
\widetilde{\rho^1_p}(J) & \text{if} & \eta=0 \text{ in } \widetilde{A_0^p}(J)\\
\widetilde{\rho^1_p}(J)+\rho^0(K) & \text{if} & \eta\neq0 \text{ in } \widetilde{A_0^p}(J)\\
\end{array}
\right.$
\end{enumerate}
\end{proposition}

\section{Strongly $p$-anisotropic knots}\label{anisotropy examples} 

For a knot $K$ the classical rational Blanchfield form $Bl$ is sesquilinear with respect to the involution $\overline{q}(t)=q(t^{-1})$.  For a symmetric polynomial $p$, this involution extends over $\widetilde{R_p}$, so the Blanchfield form extends to a sesquilinear form which we call the \textbf{strongly localized Blanchfield form},
\begin{equation*}\widetilde{Bl}_p:\widetilde{A_0^p}(K)\times \widetilde{A_0^p}(K)\to \frac{\Q(t)}{\widetilde{R_p}}.
\end{equation*}

A submodule $P$ of $\widetilde{A_0^p}(K)$ is called \textbf{isotropic} if $P\subseteq P^\perp$ with respect to $\widetilde{Bl}_p$ and it is called \textbf{Lagrangian} or \textbf{self-annihilating} if $P=P^\perp$.  A knot $K$ is called \textbf{strongly $p$-anisotropic} if $\widetilde{A_0^p}(K)$ has no nontrivial isotropic submodules.

We now provide examples of strongly $p$-anisotropic knots.  

\begin{proposition}\label{strong anisotropy}
Let $\Delta$ be the Alexander polynomial of a knot $K$.  Let $p$ be a symmetric polynomial and $h=(\Delta, p)$ be the greatest common divisor of $\Delta$ and $p$.  Suppose that $h$ has no non-symmetric factors and no roots of multiplicity greater than $1$.   If 
$\widetilde{\left(\frac{\Delta}{h},p\right)}=1$
then $K$ is strongly $p$-anisotropic.
\end{proposition}

\begin{proof}

As a first step, we show that $\widetilde{A_0^p}(K)$ is cyclic.  Since the unlocalized Alexander module $A_0(K)$ is torsion over the PID $\Q[t^{\pm1}]$, it has a decomposition into elementary factors:
\begin{equation*}
A_0(K) = \iterate{\oplus}{i=1}{k}\frac{\Q[t^{\pm1}]}{(q_i)}
\end{equation*}
where $q_i$ divides $q_{i+1}$ for each $0<i<k$ and $\displaystyle \prod_{i=1}^{k} q_i=\Delta$.  Thus, the localized Alexander module has the decomposition
\begin{equation*}
\widetilde{A_0^p}(K) =  A_0(K)\underset{\Q[t^{\pm1}]}{\otimes} \widetilde{R_p} = \iterate{\oplus}{i=1}{k}\frac{\widetilde{R_p}}{(q_i)}
\end{equation*}
If for some $i<k$ there exist some $z\in\cplx, m,n\in \Z$ such that $q_i(z^n)=p(z^m)=0$ then $q_{i+1}(z^n)$ is also zero since $q_i$ divides $q_{i+1}$.  Thus, $\Delta = \displaystyle \prod_{i=1}^{k} q_i$ has a root of multiplicity at least two at $z^n$.  Since $h$ has no roots of multiplicity greater than $1$, $\dfrac{\Delta}{h}(z^n)=0$ contradicting that $\widetilde{(\frac{\Delta}{h},p)}=1$. Thus, $q_i\in \widetilde{S_p}$ is a unit in $\widetilde{R_p}$ for every $i<k$ and 
$
\widetilde{A_0^p}(K) = \dfrac{\widetilde{R_p}}{(q_k)}.
$
In particular, $\widetilde{A_0^p}(K)$ is cyclic.

Let $\hat h = (q_k,h)$.  Since $\dfrac{q_k}{\hat h}$ divides $\dfrac{\Delta}{h}$, which is a polynomial strongly coprime to $p$, it follows that $\dfrac{q_k}{\hat h} \in \widetilde{S_p}$ is a unit in $\widetilde{R_p}$ so that the ideals generated by $\hat h$ and $q_k$ in $\widetilde{R_p}$ are the same and
$
\widetilde{A_0^p}(K) \cong \widetilde{R_p}/{(\hat h)}.
$

Let $\eta$ be a generator of $\widetilde{A_0^p}(K)$.  If $\widetilde{Bl_p}(\eta, \frac{x}{y}\eta)=0$ for some $\frac{x}{y}\in \widetilde{R_p}$, then for any other $\frac{r}{q}\eta\in\widetilde{A_0^p}(K)$, $\widetilde{Bl}_p(\frac{r}{q}\eta, \frac{x}{y}\eta)=\frac{{r}}{{q}}\widetilde{Bl}_p(\eta, \frac{x}{y}\eta)=0$.   Thus,  $\widetilde{Bl}_p(-,\frac{x}{y}\eta)$ is identically zero.  By the non-singularity of the Blanchfield form, this implies that $\frac{x}{y}\eta$ is zero, so that $\frac{x}{y}$ is zero in $\widetilde{R_p}/(\hat h)$ and $\frac{x}{y}\in (\hat{h})$

If $\frac{x}{y}\eta$ is an isotropic element of $\widetilde{A_0^p}(K)$, then 
\begin{equation*}
0=\widetilde{Bl}_p\left(\frac{x}{y}\eta, \frac{x}{y}\eta\right) = \widetilde{Bl}_p\left(\eta, \frac{x\overline x}{y \overline y}\eta\right)
\end{equation*}
so that $\frac{x\overline x}{y\overline y}\in{(\hat h)}$.  Thus, there is some $\frac{u}{v}\in \widetilde{R_p}$ such that
$
\frac{x\overline x}{y\overline y} = \hat h\frac{u}{v}
$
in $\widetilde {R_p}$.  Cross-multiplying gives the equality in $\Q[t^{\pm1}]$
\begin{equation}\label{isotropy implies}
x\overline x v = \hat h u y\overline y.
\end{equation}
Thus, $\hat h$ divides $x\overline x v$.  The fact that $v\in \widetilde{S_p}$ implies $\widetilde{(v,p)}=1$ and in particular $(v,p)=1$.  Since $\hat h$ divides $p$ it follows that $(v,\hat h)$=1.  Therefore (\ref{isotropy implies}) implies that $\hat h$ divides $x\overline x$.  Since $\hat h$ divides $h$, which has neither any non-symmetric factors nor any repeated factors, it must be that $\hat h$ has neither any non-symmetric factors nor any repeated factors.  Thus, $\hat h$ dividing $x\overline x$ implies that $\hat h$ divides $x$ so that $\frac{x}{y}$ is in the ideal of $\widetilde{R_p}$ generated by $\hat h$ and $\frac{x}{y}\eta =0$.  

Thus, $\widetilde{A_0^p}(K)$ has no nonzero isotropic submodules with respect to $\widetilde{Bl}_p$.
\end{proof}

We restrict Proposition~\ref{strong anisotropy} to the setting from which we draw examples:

\begin{corollary}\label{strong anisotropy corollary}
Let $\Delta$ be the Alexander polynomial of a knot $K$.
\begin{enumerate}
\item If $\widetilde{(p,\Delta)}=1$ then $\widetilde{A_0^p}(K)=0$ and $K$ is strongly $p$-anisotropic.
\item If $p=\Delta$ has no repeated roots and has no non-symmetric factors then $K$ is strongly $p$-anisotropic. 
\end{enumerate} 
\end{corollary}

\section{
Isotropy and extension of coefficient systems}\label{isotropy}

In this section we discuss the affect that the conditions of strong $p$-anisotropy and double anisotropy have on higher order Alexander modules of a knot.  More precisely, if a knot is strongly $p$-anisotropic we discover a restriction on the structure of isotropic submodules of certain localizations of higher order Alexander modules.  We go on to show that doubly anisotropic elements of the Alexander module produce elements of unlocalized higher order Alexander modules which are almost doubly anisotropic.


We begin by describing what we mean by higher order Alexander modules.

For a knot, $K$.  Let $\psi:\pi_1(M(K))\to \Gamma$ be a homomorphism to a PTFA group, $\Gamma$, which factors as
$$
\pi_1(M(K))\to \langle t\rangle \into A \normalsubgroup \Gamma
$$
where $t$ is the generator of the Abelianization of $\pi_1(M(K))$ and $A$ is a TFA normal subgroup of $\Gamma$.  Let 
\begin{equation*}S_p(A)=\{q_1(a_1)\dots q_n(a_n)|\widetilde{(q_i,p)}=1, a_i\in A\}\subseteq\Q[A].\end{equation*}
Since $A$ is normal, $S_p(A)$ is a $\Gamma$-invariant divisor set for $\Q[A]$, \cite[Proposition 4.1]{primaryDecomposition} shows that $S_p(A)\subseteq \Q[\Gamma]$ satisfies the right (as well as the left) Ore condition and the localization $\Rn:=\Q[\Gamma]S_p(A)^{-1}$ is defined.  For the definition of the Ore condition and a treatment of localization for noncommutative rings, see \cite[Chapter 2]{ste}.  Let $\K(\Gamma) = \Q[\Gamma](\Q[\Gamma]-\{0\})^{-1}$ be the skew field of fractions of $\Q[\Gamma]$.

We are interested in the localized higher order Alexander module of $K$, $H_1(M(K);\Rn)$.


According to \cite[Theorem 2.13]{whitneytowers}, there exists a sesquilinear form  $$Bl_\Gamma:H_1(M(K);\Rn)\times H_1(M(K);\Rn)\to\frac{\K(\Gamma)}{\Rn}.$$
We summarize the construction.  Consider the Bockstien exact sequence on cohomology,
\begin{equation*}
\begin{array}{c}
H^1(M(K);\K(\Gamma))\to H^1\left(M(K);\frac{\K(\Gamma)}{\Rn}\right)\overset{Bo}{\to} \\H^2(M(K);\Rn)\to H^2(M(K);\K(\Gamma)).
\end{array}
\end{equation*}
By Poincar\'e duality and \cite[Lemma 3.9]{C}, 
$$H^2(M(K);\K(\Gamma))\cong H_1(M(K);\K(\Gamma))=0$$
By \cite[Remark 2.8.1]{whitneytowers} there is a universal coefficient theorem for skew field coefficients and 
$$H^1(M(K);\K(\Gamma)) \cong \hom_{\K(\Gamma)}(H_1(M(K);\K(\Gamma)), \K(\Gamma)) = 0.$$
  What remains of the Bockstien exact sequence is that the Bockstien homomorphism $Bo:H^1(M(K);\K(\Gamma)/\Rn)\to H^2(M(K);\Rn)$ is an isomorphism.  The Blanchfield form, $Bl_\Gamma$, is defined by the composition
$$
\begin{array}{c}H_1(M(K);\Rn)\overset{P.D.}{\to}H^1(M(K);\Rn)\overset{Bo^{-1}}{\to}H^1(M(K);\K(\Gamma)/\Rn)\\\overset{\kappa}{\to}\hom_{\Rn}(H_1(M(K);\Rn),\K(\Gamma)/\Rn)\end{array}$$
where $P.D.$ denotes Poincar\'e duality and $\kappa$ is the Kronecker map, that is, 
$$Bl_\Gamma(a,b)=\left((\kappa\circ Bo^{-1}\circ P.D.)(a)\right)(b)$$

By \cite[Lemma 3.2 and Proposition 3.6]{Lei3}, since $\pi_1(M(K))\to \Gamma$ factors nontrivially through Abelianization and has image in the normal TFA subgroup $A$,  
$$H_1(M(K);\Rn)\cong H_1(M;\widetilde{R_p})\underset{\widetilde{R_p}}{\otimes}\Rn$$
 and for any $a\otimes \alpha$ and $b\otimes \beta$ in $H_1(M(K);\Q[\Gamma]S_p(A)^{-1})$, 
$$Bl_\Gamma(a\otimes \alpha,b\otimes \beta)=\overline{\alpha}\Psi(\widetilde{Bl_p}(a,b))\beta,$$
where $\Psi:\frac{\Q(t)}{\widetilde{R_p}}\to \frac{\K(\Gamma)}{\Rn}$ is the map induced by $\psi$.

Thus, in this section we start with a torsion $\widetilde{R_p}$ module, $M$ with a bilinear form $B:M\times M\to \frac{\Q(t)}{\widetilde{R_p}}$ and study the bilinear form on $M_\Gamma=M\underset{\widetilde{R_p}}{\otimes}\Rn$, $B_\Gamma: M_\Gamma\times M_\Gamma\to \frac{\K(\Gamma)}{\Rn}$ given by $B_\Gamma(a\otimes \alpha,b\otimes \beta)=\overline{\alpha}\Psi(B(a,b))\beta.$

\subsection{The inheritance of anisotropy under extension of coefficients.}

The following theorem reveals an aspect of the behavior of isotropic submodules under this extension of coefficients.  

\begin{theorem}\label{extension and isotropy}
Consider the infinite cyclic group $\langle t\rangle$ and a torsion $\widetilde{R_p}$ module $M$, with bilinear form $B:M\times M\to \dfrac{\Q(t)}{\widetilde{R_p}}$.  

Suppose $\langle t\rangle$ injects into a TFA group $A$ which is a normal subgroup of a PTFA group $\Gamma$.  If $P$ is an isotropic submodule of $M\otimes \Rn$ with respect to $B_\Gamma$, then $\{m\in M | m\otimes 1\in P\}$ is isotropic with respect to $B$.
\end{theorem}

Theorem \ref{extension and isotropy} is a consequence of Lemma~\ref{isotropy extension} below.

\begin{lemma}\label{isotropy extension}
Suppose $t\mapsto T\in A$ defines a monomorphism from $\langle t\rangle$ to $A$ where $A$ is a TFA group and a normal subgroup of a PTFA group $\Gamma$.  Then the induced map
\begin{equation*}\displaystyle\Psi:\dfrac{\Q(t)}{\widetilde{R_p}} \into \dfrac{\K(\Gamma)}{\Rn}.\end{equation*}
is a monomorphism.
\end{lemma}
\begin{proof}

Suppose that $\displaystyle\dfrac{f(t)}{g(t)}$ is in the kernel of this map.  Then $\displaystyle\dfrac{f(T)}{g(T)}$ is contained in $\Rn$ and there exist some $r\in \Q[\Gamma]$ and $q\in S_p(A)$ such that $\displaystyle\dfrac{f(T)}{g(T)} = \dfrac{r}{q}$.  By the definition of equality in $\K(\Gamma)$ (see \cite[Chapter 2 Proposition 1.4]{ste}) This implies that there exist nonzero $c,d\in \Q[\Gamma]$ such that 
\begin{eqnarray}
\label{E1}f(T)c&=&rd,\\
\label{E2}g(T)c&=&qd.
\end{eqnarray}

Considering (\ref{E1}) as an equation in $\Q[\Gamma](\Q[A]-\{0\})^{-1}$ (into which $\Q[\Gamma]$ injects), it reduces to $c=(f(T))^{-1}rd$.  This substitution reduces (\ref{E2}) to $g(T)(f(T))^{-1}rd=qd$.  Cancelation gives that $g(T)(f(T))^{-1}r=q$.  Now, $g(T)$ and $(f(T))^{-1}$ sit in the image of the field $\Q(t)$ and so commute with each other.  Thus, $(f(T))^{-1}g(T)r=q$.  Multiplying by $f(T)$ on the left gives us that 
\begin{equation}\label{E2'}g(T)r=f(T)q.\end{equation}

Let $X$ be a transversal for $\Gamma/A$ (that is, a subset of $\Gamma$ containing the identity, $1$, such that every equivalence class in $\Gamma/A$ has a unique representative in $X$).  The group ring $\Q[\Gamma]$ is free as a $\Q[A]$ module and has basis $X$, so that $r$ can be uniquely realized as $r=\displaystyle \sum_{x\in X}r_x x$ where each $r_x$ is in $\Q[A]$.  Since the right hand side of (\ref{E2'}) is in $\Q[A]$, that is, the span of $\{1\}$, it reduces to
\begin{eqnarray}
\label{E2''}
g(T)r_1=f(T)q\text{ and }\\
r_x=0\text{ for }x\in X-\{1\}.
\end{eqnarray}

Now suppose that $\displaystyle \dfrac{f(t)}{g(t)} \in \Q(t)$ is in reduced terms and that $g$ in not in $\widetilde{S_p}$.  Since $q$ is assumed to be in $S_p(A)$, there exist polynomials $q_1,\dots, q_n\in \widetilde{S_p}$ and $a_1,\dots, a_n\in A$ so that $q=\displaystyle \prod_{i=1}^n q_i(a_i)$. Equation (\ref{E2''}) as an equality in $\Q[A]$ involves only finitely many elements of $A$, namely, $T, a_1,\dots, a_n$ and $b_1,\dots, b_k$ where $b_1,\dots ,b_k$ are whatever terms appear in $r_1$.  The span of $\{T,a_1,\dots, a_n, b_1,\dots, b_k\}$ is a finitely generated subgroup of the TFA group, $A$, so it is free Abelian.  Pick a basis $\{s, c_1,\dots, c_m\}$ such that $s^l=T$ for some $l$.  The equality (\ref{E2''}) can then be realized as an equality in the multivariable Laurent polynomial ring $\Q[s^{\pm1}, c_1^{\pm1},\dots, c_m^{\pm1}]$: 
\begin{equation}\label{E23}
g(s^k)r_1(s,c_1,\dots, c_m) = f(s^k)\prod_{i=1}^n q_i(s^{k_i} c_1^{k_{i,1}} \dots c_m^{k_{i,m}}).
\end{equation}

Since $g$ is not an element of $\widetilde{S_p}$, (i.e. it is not strongly coprime to $p$) there exists some $z\in \cplx$ and $\alpha,\beta\in \Z-\{0\}$ such that $g(z^\beta)=p(z^\alpha)=0$.  Evaluating (\ref{E23}) at $s=z^{\beta/k}$ and $c_1=\dots=c_n=1$ gives an equality in $\cplx$
\begin{equation}\label{evaluated}
0 = f(z^\beta)\prod_{i=1}^n q_i(z^{\beta k_i/k}).
\end{equation}
Since $z^\beta$ is a root of $g$ and $f$ is assumed to be relatively prime to $g$, $f(z^\beta)\neq 0$.  Therefore $q_i(z^{\beta k_i/k})=0$  for some $i$, contradicting that $q_i\in \widetilde{S_p}$ is strongly coprime to $p$.  Thus, it must be that $g\in \widetilde{S_p}$ and $\frac{f}{g}=0$ in $\Q(t)/\widetilde{R_p}$.
\end{proof}

We now prove Theorem~\ref{extension and isotropy}.

\begin{proof}[Proof of Theorem~\ref{extension and isotropy}]

Consider any $m,n\in M$ such that $m\otimes1,n\otimes 1\in P$.  Since $P$ is isotropic,  $0=B_\Gamma(m\otimes1,n\otimes 1)=\Psi(B(m,n))$.  Since $\Psi$ is injective by Lemma~\ref{isotropy extension} this implies that $B(m,n)=0$ completing the proof.
\end{proof}

\subsection{The inheritance of double anisotropy under extension of coefficients}

The following Proposition reveals that the double anisotropy is partially inherited under the extension of coefficients.  A similar claim could be proven after localization.  For our purposes the unlocalized claim is sufficient.

\begin{proposition}\label{2xAnisoProp}
Let $A$ be a TFA group which is a normal subgroup of a PTFA group $\Gamma$.  Let $M$ be a torsion $\Q[A]$ module with $\Q[A]$-sesquilinear form $B_A:M\times M\to \dfrac{\K(A)}{\Q[A]}.$
  This sesquilinear form extends to a $\Q[\Gamma]$-sesquilinear form $B_\Gamma:M\otimes \Q[\Gamma]\times M\otimes \Q[\Gamma]\to \dfrac{\K(\Gamma)}{\Q[\Gamma]}.$
  
    Let $Q$ be an isotropic submodule of $M$ with respect to $B_A$ and $P$ be an isotropic submodule of $M\otimes \Q[\Gamma]$ with respect to $B_\Gamma$.  If $\eta\otimes 1\in M\otimes \Q[\Gamma]$ sits in $P + Q\otimes \Q[\Gamma]$, then $\eta=p+q$ for some $q\in Q$ and $B_A(p,p)=0$.
\end{proposition}
\begin{proof}
Let $X$ be a transversal for $\Gamma/A$.  Suppose that \begin{equation}\label{isotropy assumption}
\eta\otimes 1 = \sum_{x\in X}p_x \otimes x + \sum_{x\in X}q_x \otimes x 
\end{equation}
with $q_x\in Q$, $p_x\in M$ for all $x$ and $p:=\displaystyle  \sum_{x\in X}p_x \otimes x \in P$.  Thinking of $\Q[\Gamma]$ as the free $\Q[A]$ module generated by $X$, (\ref{isotropy assumption}) implies that $p_x+q_x=0$ for $x\neq 1$ and $\eta=p_1+q_1$.  Since $P$ is isotropic, 
\begin{equation}\label{p isotropic 1}
\displaystyle
0=B_\Gamma(p,p)=B_\Gamma\left(\sum_{x\in X} p_x\otimes x, \sum_{y\in X} p_y\otimes y\right).
\end{equation}
Appealing to the $\Gamma$-sesquilinearity of $B_\Gamma$, this implies 
\begin{eqnarray}\label{p isotropic 2}
\begin{array}{rcl}
0&=&\sum_{x\in X} \sum_{y\in X} x^{-1} B_\Gamma( p_x\otimes 1, p_y\otimes1)y
\\
&=&\sum_{x\in X} \sum_{y\in X} (x^{-1} B_A( p_x, p_y)x)x^{-1}y.
\end{array}
\end{eqnarray}
Since $A$ is normal in $\Gamma$, $\displaystyle x^{-1} B_A( p_x, p_y)x$ is in  $\frac{\K(A)}{\Q[A]}$ for each $x, y\in X$.

Since $X$ is a choice of coset representatives for $\Gamma/A$, each $x^{-1}y$ is equivalent modulo $A$ to some $z$ in $X$, that is, there is some $a_{x,y}\in A$ such that $x^{-1}y=a_{x,y}z$.  We use this to rearrange (\ref{p isotropic 2}), 
\begin{equation}\label{p isotropic 3}
\displaystyle
0=\sum_{z\in X} \sum_{x^{-1}y\equiv z} (x^{-1} B_A( p_x, p_y)x)a_{x,y}z
\end{equation}
The map $\left(\dfrac{\K(A)}{\Q[A]}\right)^{|X|}\to \dfrac{\K(\Gamma)}{\Q[\Gamma]}$ defined by sending $\langle r_x\rangle_{x\in X}$ to $\displaystyle\sum_{x\in X}r_x\otimes x$ is injective.  Indeed, if $\langle \frac{a_x}{b_x}\rangle_{x\in X}$ is in the kernel of this map then there exists some $c=\displaystyle\sum_{x\in X}c_x x$ with $c_x\in \Q[A]$ such that 
$$\displaystyle\sum_{x\in X}a_x b_x^{-1} x=\displaystyle\sum_{x\in X}c_x x$$
Since $A$ is Abelian, $a_x b_x^{-1}=b_x^{-1} a_x$.  Left multiplying by $b=\displaystyle\prod_{x\in X}b_x$ we see that
$$\displaystyle\sum_{x\in X}(b_x^{-1}b)a_x x=\displaystyle\sum_{x\in X}b c_x x$$
This is an equation in $\Q[\Gamma]$.  The set $X$ is a basis for $\Q[\Gamma]$ as a free $\Q[A]$ module.  Thus, for all $x\in X$,
$(b_x^{-1}b)a_x=b c_x$ so $\frac{a_x}{b_x}=c_x\in \Q[A]$ holds in $\Q(A)$ and $\langle \frac{a_x}{b_x}\rangle_{x\in X}$ is zero in $\left(\dfrac{\K(A)}{\Q[A]}\right)^{|X|}$.

This injectivity together with (\ref{p isotropic 3}) implies that for each $z\in X$, 
\begin{equation}\label{p isotropic 4}
\displaystyle
0=\sum_{x^{-1}y\equiv z} (x^{-1} B_A( p_x, p_y)x)a_{x,y}.
\end{equation}
Taking $z=1\in X$ we see that 
\begin{equation}\label{p isotropic 5}
\displaystyle
0=\sum_{x\in X} (x^{-1} B_A( p_x, p_x)x)a_{x,x}.
\end{equation}
As we observed previously, for $x\neq 1$, $p_x=-q_x$ is in $Q$, so $B_A(p_x, p_x)=0$.  Thus, all but one of the terms in (\ref{p isotropic 5}) vanishes.  Dropping them, we see that 
$
0=B_A(p_1, p_1).
$
Since, $\eta=p_1+q_1$, this completes the proof.
\end{proof}

\begin{proposition}\label{2xAnisoCor}
Let $M$ be a torsion $\Q[t^{\pm1}]$-module with bilinear form $B:M\times M\to \frac{\Q(t)}{\Q[t^{\pm1}]}$.  Let $Q\subseteq M$ be isotropic.

Suppose $\langle t\rangle$ injects into a TFA group $A$ which is  a normal subgroup of a PTFA group $\Gamma$ and $P\subseteq M\otimes \Q[\Gamma]$ is isotropic.  If $\eta\otimes 1$ is in $P+ Q\otimes \Q[\Gamma]$, then $\eta$ is not doubly anisotropic with respect to $B$.
\end{proposition}
\begin{proof}
Consider $\eta^A=\eta\otimes 1\in M\otimes \Q[A]$ and the isotropic module $Q^A=Q\otimes \Q[A]$.  Notice that $\eta^A\otimes 1\in (M\otimes \Q[A])\otimes \Q[\Gamma]$ is in $Q^A\otimes \Q[\Gamma]+P$ by assumption.  Applying Proposition~\ref{2xAnisoProp} gives that there is some $q^A \in Q^A$ and $p^A\in M\otimes \Q[A]$ with $B_A(p^A,p^A)=0$ and $\eta^A=p^A+q^A$.

But then $\eta\otimes 1\in M\otimes \Q[A]$  sits in the sum of $Q\otimes \Q[A]$ with the isotropic submodule $\langle p^A\rangle$.  Applying Proposition~\ref{2xAnisoProp} again gives that $\eta=p+q$ where $B(p,p)=0$ and $q\in Q$.  In this case, $\eta$ sits in the sum of the isotropic submodules $Q$ and $\langle p\rangle$ and so is not doubly anisotropic.
\end{proof}


\section{The proof of Theorem~\ref{big corollary}}\label{main section}

In this section we set out to prove the main result of this paper, Theorem~\ref{big corollary}:

\begin{reptheorem}{big corollary}
Let $\{K_i\}$ be a possibly infinite set of knots:
\begin{enumerate}
\item whose Alexander polynomials are strongly coprime,
\item whose Tristram-Levine signatures have vanishing integrals,
\item whose prime factors have square-free Alexander polynomials and
\item whose $\rho^1$-invariants do not vanish, that is $\rho^1(K_i)\neq 0$ .
\end{enumerate}
   For $i=1,2, \dots$ and $j=1,2,\dots, n$ let $R_{i,j}$ be a slice knot and $\eta_{i,j}$ be an unknotted curve in the complement of $R_{i,j}$ such that the pair $(R_{i,j},\eta_{i,j})$ is doubly anisotropic.

Let $K_i^0=K_i$ and $K_i^j=R_{i,j}(\eta_{i,j},K_i^{j-1})$.

 Then $\left\{K_i^n\right\}_{i=1}^\infty$ is linearly independent modulo ($n+1.5$)-solvable knots.
\end{reptheorem}

In order to prove the theorem we explore the interaction between an iterated infection procedure, strongly localized $\rho$-invariants and the following technical condition on a bounded 4-manifold.

\begin{definition}
Let $K_1, \dots, K_m$ be knots in $S^3$.  Consider a 4-manifold $W$ with $\bdry W = \sqcup M(K_i)$ and an epimorphism $\phi:\pi_1(W)\onto\Gamma$.  The pair $(W,\Gamma)$ is said to satisfy \textbf{condition C} with respect to integers $n, h$ (which we abbreviate by  saying that $(W,\Gamma)$ is $\mathbf{C(n,h)}$) if the following conditions hold:
\begin{enumerate}[label= ({C}\arabic*)]
\item\label{C1} $\Gamma^{(n+1)}_\mathfrak{r}=0$ (This condition implies that $\Gamma$ is PTFA).
\item\label{C2} There is a normal Abelian subgroup $A\lhd\Gamma$ such that for each $i$, there is a monomorphism $\alpha_i:\nquotient{M(K_i)}{1}\cong \Z \to\Gamma$ making the following diagram commute,
\begin{equation*}
	\begin{diagram}
	\node{\pi_1(M(K_i))} \arrow{e}\arrow{s}\node{\nquotient{M(K_i)}{1}}\arrow{s,r,J}{\alpha_i}\\
	\node{\pi_1(W)} \arrow{e,t}{\phi}\node{\Gamma}
	\end{diagram},
\end{equation*}
and $\im(\alpha_i)$ sits in $A$.  The subgroup $A$ does not depend on $i$.
\item\label{C3} For any coefficient system $\psi:\pi_1(W)\to\Lambda$ with $\Lambda^{(h+1)}_\mathfrak{r}=1$ and an epimorphism, $\beta$, making the following diagram commute

\begin{equation*}
	\begin{diagram}
	\node{\pi_1(W)} \arrow{e,t}{\phi}\arrow{se,l}{\psi}
	\node{\Gamma} \\
	\node[2]{\Lambda}\arrow{n,r,A}{\beta}
	\end{diagram}
\end{equation*}
and for any Ore localization $\Q[\Lambda]S^{-1}$ of $\Q[\Lambda]$, \begin{equation*}\ker\left(H_1(\bdry W;\Q[\Lambda]S^{-1})\to H_1(W;\Q[\Lambda]S^{-1})\right)\end{equation*} is isotropic with respect to the Blanchfield form $$Bl_\Lambda:H_1(\bdry W;\Q[\Lambda]S^{-1})\times H_1(\bdry W;\Q[\Lambda]S^{-1})\to \frac{\K(\Lambda)}{\Q[\Lambda]S^{-1}}.$$
\item\label{C4} For any PTFA coefficient system $\psi:\pi_1(W)\to\Theta$ on $W$, with $\Theta_\mathfrak{r}^{(h+2)}=0$, $\sigma^{(2)}(W;\psi)-\sigma(W)=0$
\end{enumerate}

If such a pair $(W,\Gamma)$ exists then we say that $\sqcup M(K_i)$ bounds a $C(n,h)$.
\end{definition}

\begin{remark}\label{C3 to isotropy}
Notice that if $n\le h$, then taking $\Lambda=\Gamma$, $\beta$ to be the identity map and $A$ to be the Abelian group of condition \ref{C2}, then condition $\ref{C3}$ implies that for a polynomial, $p$, \begin{equation*}\ker\left(H_1(\bdry W;\Q[\Gamma]S_p(A)^{-1})\to H_1(W;\Q[\Gamma]S_p(A)^{-1})\right)\end{equation*}
is isotropic with respect to the Blanchfield form on $H_1(\bdry W;\Q[\Gamma]S_p(A)^{-1})$.  
\end{remark}

Notice that ($h.5$)-solutions satisfy this condition.

\begin{lemma}\label{base case}
If $W$ is an $h.5$ solution for $K$ and $\phi:\pi_1(W)\to \Z$ is the Abelianization map then $(W,\Z)$ is $C\left(0, h-1\right)$.
\end{lemma}
\begin{proof}
Condition \ref{C1} holds since $\Z$ is torsion free Abelian. Condition \ref{C2} holds since by definition of $h.5$ solvability \cite[Definition 1.2]{whitneytowers} the inclusion induced map on first homology is an isomorphism.  Condition \ref{C3} follows from \cite[Theorem 6.3]{blanchfieldDuality} since $W$ being an $(h.5)$-solution implies that it is an $(h)$-solution and so is an $(h)$-null-bordism.  Condition \ref{C4} follows from \cite[Theorem 4.2]{whitneytowers}.  
\end{proof}

The action of connected sum and infection on Condition $C$ are provided by the following lemmas.  Their proofs are delayed until subsection~\ref{proofs}.

\begin{lemma}\label{connected sum inheritance}
Let $K_i=K_{i,1} \# \dots \#K_{i,m_i}$ for $i=1, \dots, p$.  If $\iterate{\sqcup}{i=1}{p} M(K_i)$ bounds a $C(n,h)$, $(W,\Gamma)$, then $\iterate{\sqcup}{i=1}{p}\iterate{\sqcup}{j=1}{m_i}M(K_{i,j})$ bounds a $C(n,h)$.\end{lemma}

\begin{lemma}\label{infection inheritance}
  For $1\le i\le m$, let $K_i$ be a knot, $R_i$ be a slice knot and $\eta_i$ be an unknotted curve representing a doubly anisotropic element of $A_0(R_i)$.  If $\iterate{\sqcup}{i=1}{m}M(R_i(\eta_i,K_i))$ bounds a $C(n,h)$ with $n \le h$, then $\iterate{\sqcup}{i=1}{m}M(K_i)$ bounds a  $C(n+1,h)$.
\end{lemma}

We now combine the three lemmas above to discover a relationship between the result of iterated infection being solvable and zero surgery on the deepest infecting knots cobounding a $4$-manifold satisfying condition $C$.

\begin{lemma}\label{solvable to C}
For integers $1\le i\le m$ and $1\le j\le n$ let $R^{i,j}$ be a slice knot with doubly anisotropic curve $\eta_{i,j}$.  For $1\le i\le m$ let $K_i$ be a knot.  Let $K^{i,j}$ be recursively defined by $K^{i,0} = K_i$ and $K^{i,j} = R^{i,j}_{\eta_{i,j}}(K^{i,j-1})$.  If $\iterate{\#}{i=1}{m}K^{i,n}$ is ($h.5$)-solvable for $h\ge n$, then $\iterate{\sqcup}{i=1}{m}M(K_{i})$ bounds a 4-manifold with coefficient system $\Gamma$ such that $(W,\Gamma)$ is $C(n,h-1)$.  

\end{lemma}
\begin{proof}

If $\iterate{\#}{i=1}{m}K^{i,n}$ is ($h.5$)-solvable then Lemma~\ref{base case} implies that $M\left(\iterate{\#}{i=1}{m}K^{i,n}\right)$ bounds a $C(0,h-1)$.  Applying Lemma~\ref{connected sum inheritance}, this means that $\iterate{\sqcup}{i=1}{m} M\left(K^{i,n}\right)$ bounds a $C(0,h-1)$. 

Now, applying Lemma~\ref{infection inheritance} (if $h-1\ge 0$) implies that $\iterate{\sqcup}{i=1}{m}  M\left(K^{i,n-1}\right)$ bounds a $C(1,h-1)$.  Applying it again (if $h-1\ge 1$) implies that $\iterate{\sqcup}{i=1}{m}  M\left(K^{i,n-2}\right)$ bounds a $C(2,h-1)$.  Applying it a total of $n$ times (provided that $h-1\ge n-1$) gives that $\iterate{\sqcup}{i=1}{m} M\left(K^{i,n-n}\right) = \iterate{\sqcup}{i=1}{m} M\left(K_{i}\right)$ bounds a $C(n,h-1)$, as claimed.
\end{proof}

Next, we prove that the $\widetilde{\rho^1_p}$-invariant is an obstruction to a disjoint union of zero surgeries on knots bounding a $4$-manifold satisfying condition $C$.

\begin{lemma}\label{C to rho}
Let $p\in \Q[t^{\pm1}]$ be a polynomial.  If $\{K_i\}_{i=1}^m$ is a set of knots such that for each $i$, $K_i$ decomposes as a connected sum of $p$-anisotropic knots and $\iterate{\sqcup}{i=1}{m}M(K_i)$ bounds a 4-manifold $W$ with coefficient system $\Gamma$ such that $(W,\Gamma)$ is $C(n,h)$ with $n\le h$, then $\displaystyle\sum_{i=1}^{m}\widetilde{\rho^1_p}(K_i)=0$.
\end{lemma}

\begin{proof}
We first prove the lemma in the more restrictive setting that each $K_i$ is strongly $p$-anisotropic.  By \ref{C2}, the coefficient system $\phi:\pi_1(W)\to \Gamma$ restricted to the $M(K_i)$-boundary component factors nontrivially through the Abelianization of $\pi_1(M(K_i))$ and by \ref{C4} the associated signature defect is zero.  Thus,
\begin{equation*}
\displaystyle \sum_i \rho^0(K_i)=0.
\end{equation*}  
While noteworthy, this is not the desired conclusion.  We find another coefficient system on $W$ which, when restricted to each $M(K_i)$ boundary component, factors injectively through the quotient of $\pi_1(M(K_i))$ by $\pi_1(M(K_i))^{(2)}_{\widetilde p}$.  

By Remark~\ref{C3 to isotropy} \begin{equation*}P=\ker(H_1(M(K_i);\Q[\Gamma]S_p(A)^{-1})\to H_1(W;\Q[\Gamma]S_p(A)^{-1}))\end{equation*}
is isotropic.  

Since $\pi_1(M(K))\to \Gamma$ factors nontrivially through the Abelianization, \cite[Lemma 3.2 and Proposition 3.6]{Lei3} reveals that
$$H_1(M(K_i);\Q[\Gamma]S_p(A)^{-1})\cong H_1(M(K_i);\widetilde{R_p})\underset{\widetilde{R_p}}{\otimes}\Q[\Gamma]S_p(A)^{-1},$$
and for any $a\otimes \alpha$ and $b\otimes \beta$ in $H_1(M(K_i);\widetilde{R_p})\underset{\widetilde{R_p}}{\otimes}\Q[\Gamma]S_p(A)^{-1}$, $$Bl_\Gamma(a\otimes \alpha, b\otimes \beta)=\overline{\alpha} \Psi(\widetilde{Bl_p}(a,b))\beta,$$ where $\Psi:\frac{\Q(t)}{\widetilde{R_p}}\to \frac{\K(\Gamma)}{\widetilde{R_p}}$ is induced by the map $\alpha_i$ of condition \ref{C2}.  Thus, we can think of $P$ as an isotropic submodule of $\widetilde{A_0^p}(K)\otimes \Q[\Gamma]S_p(A)^{-1}$ with respect to the bilinear form induced by $\widetilde{Bl_p}$

By applying Theorem~\ref{extension and isotropy}, we see that the kernel of the map
\begin{equation}\label{injective}
\begin{array}{c}
\widetilde{A_0^p}(K_i)\to \widetilde{A_0^p}(K_i)\underset{\widetilde{R_p}}{\otimes}\Q[\Gamma]S_p(A)^{-1} = H_1(W;\Q[\Gamma]S_p(A)^{-1}),
\end{array}
\end{equation} which is equal to $\left\{m\in\widetilde{A_0^p}(K_i)|m\otimes 1\in P\right\}$ is isotropic.  We assume that $\widetilde{A_0^p}(K_i)$ has no nontrivial isotropy so \eref{injective} is injective.

Now we build a new coefficient system on $W$.  Let $G:=\ker\left(\phi:\pi_1(W)\to \Gamma\right)$, so that $G$ is isomorphic to the fundamental group of $\widetilde{W}_\Gamma$, the $\Gamma$-cover of $W$.  Let $G^{(1)}_p\subseteq G$ be given by the kernel of the following composition:

\begin{equation*}
F:{G}\overset{\cong}{\to}{\pi_1(\widetilde{W_\Gamma})}\to{\dfrac{H_1(\widetilde{W_\Gamma};\Z)}{\Z\text{-torsion}}}\hookrightarrow{H_1(W;\Q[\Gamma])}\to{H_1(W;\Q[\Gamma]S_p(A)^{-1})}
\end{equation*}

Observe that $G^{(1)}_p$ is normal in $\pi_1(W)$.  In order to see this, let $g\in G^{(1)}_p$ and $\gamma\in \pi_1(W)$.  Then $F(\gamma^{-1} g \gamma)$ is given by letting $\gamma_*$, the deck translation on $\widetilde{W_\Gamma}$ corresponding to $\gamma$, act on $F(g)\in H_1(W;\Q[\Gamma]S_p(A)^{-1})$.  Since $F(g)=0$, it follows that $F(\gamma^{-1}g\gamma)=\gamma_*(F(g))=0$, so $\gamma^{-1}g\gamma$ is in $G^{(1)}_p$.

Since \ref{C1} gives us that $\pi_1(M(K_i))\to\Gamma$ factors through Abelianization, it follows that $\phi$ is trivial on $\pi_1(M(K_i))^{(1)}$ and the map induced by inclusion sends $\pi_1(M(K_i))^{(1)}$ to $G$.  Consider the following commutative diagram:
\begin{equation}\label{diagram}
\begin{diagram}
\node{\pi_1(M(K_i))^{(1)}}\arrow{e,t}{a}\arrow{s,l}{i_*}\node{\dfrac{\pi_1(M(K_i))^{(1)}}{\pi_1(M(K_i))_{\widetilde p}^{(2)}}}\arrow{e,t,J}{b}\arrow{s,r,..}{\beta}\node{\widetilde{A_0^p}(K)}\arrow{s,l,J}{c}
\\
\node{G}\arrow{e,t}{e}\node{\dfrac{G}{G^{(1)}_p}}\arrow{e,t,J}{d}\node{H_1(W;\Q[\Gamma]S_p(A)^{-1})}
\end{diagram}
\end{equation}

The dotted map, $\beta$, is induced by $i_*$.  In order to see that it is well defined, one must check that $i_*$ maps ${\pi_1(M(K_i))_{\widetilde p}^{(2)}}$ to $G^{(1)}_p$.  In order to see this, take  $x\in \pi_1(M(K_i))_{\widetilde p}^{(2)}$.  It follows that $c(b(a(x)))$ is zero in $H_1(W;\Q[\Gamma]S_p(A)^{-1})$.  By the commutativity of the diagram, $d(e(i_*(x)))$ is zero in $H_1(W;\Q[\Gamma]S_p(A)^{-1})$ and so $i_*(x)\in \ker(e)=G^{(1)}_p$.  Thus, the map $\beta$ is well defined.

The maps $b$ and $d$ in \eref{diagram} are injections by the definition of $\pi_1(M(K))^{(2)}_{\widetilde p}$ and $G_p^1$.  The map $c$ is the monomorphism in (\ref{injective}).  Thus, $\beta$ is injective

If $x$ is in the kernel of the composition 
\begin{equation}\label{coeff sys 1}
\pi_1(M(K_i))\overset{i_*}{\to}\pi_1(W)\to\dfrac{\pi_1(W)}{G^{(1)}_p},
\end{equation}
then $i_*(x)$ is in $G^{(1)}_p\subseteq G:=\ker(\pi_1(W)\to\Gamma)$.  Since the map from $\pi_1(M(K_i))$ to $\Gamma$ factors nontrivially through Abelianization, it must be that $x$ is in $\pi_1(M(K))^{(1)}$.  This means that $x\in \ker\left(\pi_1(M(K_i))^{(1)}\underset{i_*}{\to} G\to \dfrac{G}{G^{(1)}_p}\right)$.  By the commutativity of \eref{diagram} and the injectivity of $\beta$, the kernel of this map is $\pi_1(M(K_i))^{(2)}_{\widetilde p}$.  

If we set $\Theta:=\dfrac{\pi_1(W)}{G^{(1)}_p}$ then the following commutative diagram holds for each $M(K_i)$-boundary component:
\begin{equation*}
\begin{diagram}
\node{\pi_1(M(K_i))}\arrow{s}\arrow{e}\node{\dfrac{\pi_1(M(K_i))}{\pi_1(M(K_i))^{(2)}_{\widetilde p}}}\arrow{s,J}\\
\node{\pi_1(W)}\arrow{e}\node{\Theta.}
\end{diagram}
\end{equation*}
  This implies that $\displaystyle \sum_{i=1}^m \widetilde{\rho^1_p}(K) = \sigma^2(W,\Theta)-\sigma(W)$.  It remains only to check that this signature defect is zero.

In order to apply condition \ref{C4} to get this conclusion, consider the following short exact sequence:
\begin{equation*}
0\to\dfrac{G}{G^{(1)}_p}\to \dfrac{\pi_1(W)}{G^{(1)}_p}\to\dfrac{\pi_1(W)}{G}\to 0
\end{equation*}
The leftmost term, $\dfrac{G}{G^{(1)}_p}$, is TFA.  The rightmost term, $\dfrac{\pi_1(W)}{G}$, injects into $\Gamma$, so that $\left(\dfrac{\pi_1(W)}{G}\right)^{(n+1)}_{\mathfrak r}=0$ since $\Gamma^{(n+1)}_\mathfrak{r}=1$.  This implies that 
\begin{equation*}
\Theta^{(n+2)}_{\mathfrak r}=\left(\dfrac{\pi_1(W)}{G^{(1)}_p}\right)^{(n+2)}_{\mathfrak r}=0
\end{equation*}
and condition \ref{C4} applies to give that $\sigma^2(W,\Theta)-\sigma(W)=0$.  This completes the proof in the case that each $K_i$ is $p$-anisotropic.  

In order to see it under the weaker assumption that each $K_i$ has only strongly $p$-anisotropic factors, suppose that $K_i=\iterate{\#}{b=1}{B_i}J_{i,b}$ with each $J_{i,b}$ strongly $p$-anisotropic.  An application of Lemma~\ref{connected sum inheritance} gives that 
\begin{equation*}
\iterate{\sqcup}{i=1}{m}\iterate{\sqcup}{b=1}{B_i}M(J_{i,b}).
\end{equation*}
bounds a $C(n, h)$.  Now we can apply the theorem in the case already proven to see that
\begin{equation*}
\sum_{i=1}^{m}\left(\sum_{b=1}^{B_i}\widetilde{\rho^1_p}(J_{i,b})\right)=0
\end{equation*}
By Proposition~\ref{homomorphism}, $\displaystyle\sum_{b=1}^{B_i}\widetilde{\rho^1_p}(J_{i,b})= \widetilde{\rho^1_p}(\iterate{\#}{b=1}{B_i}J_{i,b}) = \widetilde{\rho^1_p}(K_i)$.  Making this substitution completes the proof.

\end{proof}

We are now ready to prove Theorem~\ref{big corollary}.  We prove a stronger theorem from which we get it as a corollary.

\begin{theorem}\label{big theorem}
Let $p$ be a polynomial.  Let $\{K_i\}$ be a possibly infinite set of knots each of which decomposes into a connected sum of strongly $p$-anisotropic knots.   For $i=1,2, \dots$ and $j=1,2,\dots, n$ let $R_{i,j}$ be a slice knot and $\eta_{i,j}$ be an unknotted curve in the complement of $R_{i,j}$ representing a doubly anisotropic element of $A_0(R_{i,j})$.  

Let $K_i^0=K_i$ and $K_i^j=R_{i,j}(\eta_{i,j},K_i^{j-1})$.

 If $\iterate{\#}{i=1}{m}a_iK_i^n$ is (${n+1.5}$)-solvable then $\displaystyle \sum_{i=1}^{m} a_i\widetilde{\rho^1_p}(K_i)=0$.

\end{theorem}
\begin{proof}
Suppose that $\iterate{\#}{i=1}{m} a_i K_i^n$ were ($n+1.5$)-solvable.  By Lemma~\ref{solvable to C}  $\iterate{\sqcup}{i=1}{m} \left(\iterate{\sqcup}{k=1}{a_i} M(K_i)\right)$ bounds a $C(n,n)$ so by Lemma~\ref{C to rho} $\displaystyle \sum_i a_i\widetilde{\rho^1_p}(K_i)=0$.  

\end{proof}

\begin{proof}[Proof of Theorem~\ref{big corollary}]
Suppose that some linear combination $\iterate{\#}{j=1}{m} a_i K_i$ is $(n+1.5)$-solvable.  Let $p$ be the Alexander polynomial of $K_i$.  For $j\neq i$ $K_j$ has Alexander polynomial strongly coprime to $p$.  By Corollary~\ref{strong anisotropy corollary}, $K_j$ is strongly $p$-anisotropic for all $j$.  Theorem~\ref{big theorem} gives that 
\begin{equation}\label{almost done}
\sum_{j}a_j\widetilde{\rho^1_p}(K_j)=0
\end{equation}
Proposition~\ref{easy} applies to give that $\widetilde{\rho^1_p}(K_i)=\rho^1(K_i)\neq 0$ and that for $j\neq i$, $\widetilde{\rho^1_p}(K_j)=\rho^0(K_j)= 0$.  plugging these into \eref{almost done} yields $a_i\rho^1(K_i)=0$ so $a_i=0$.  

Since the choice of $i$ was arbitrary, $a_i=0$ for all $i$ and there are no nontrivial linear relationships amongst these knots modulo $n+1.5$ solvability.
\end{proof}

\subsection{Proofs of Lemmas~\ref{connected sum inheritance} and \ref{infection inheritance}}\label{proofs}

Before we prove these two important lemmas we discuss the cobordisms used to prove them

\begin{definition}
For knots $K_1, \dots, K_n$, let $V_\#={V_\#(K_1, \dots, K_n)}$ be the cobordism between $\sqcup M(K_i)$ and $M(\# K_i)$ constructed by starting with $\iterate{\sqcup}{i=1}{n} M(K_i)\times[0,1]$ and connecting it by gluing together neighborhoods of curves in $M(K_{i-1})\times\{1\}$ and $M(K_i)\times\{1\}$ representing the meridians of $K_{i-1}$ and $K_i$.
\end{definition}

\begin{definition}
Consider knots $K$ and $J$ and an unknotted curve, $\eta$, in $S^3-K$ which has zero linking with $K$.  Let $V_\infect=V_\infect(K,\eta,J)$ be the cobordism between $M(K)\sqcup M(J)$ and $M(K_\eta(J))$ given by starting with $M(K)\times I\sqcup M(J)\times I$ and gluing a neighborhood of $\eta$ in $M(K)\times\{1\}$ to a neighborhood of the meridian of $J$ in $M(J)\times \{1\}$.
\end{definition}

By virtue of the fact that the inclusion induced maps $H_2(M(K))\oplus H_2(M(J))\to H_2(V_\#)$ and $H_2(M(R))\oplus H_2(M(K))\to H_2(V_\infect)$ are epimorphisms, each of $V_\#$ and $V_\infect$  are rational $(k)$-null-bordisms for every nonnegative integer, $k$ (see \cite[Definition 5.1]{blanchfieldDuality}).  

The following is a key result about rational $(k)$-null-bordisms.  (For an $\Rn$ module $M$, $T(M)$ denotes the $\Rn$-torsion part of $M$.)

\begin{theorem*}[  {\cite[Theorem 6.3]{blanchfieldDuality}}  ]
Suppose $W$ is a rational ($k$)-null-bordism and $\phi:\pi_1(W)\to\Gamma$ is a nontrivial coefficient system where $\Gamma$ is a PTFA group with $\Gamma^{(k)}=1$.  Let $\Rn$ be an Ore localization of $\Z[\Gamma]$ so $\Z[\Gamma]\subseteq \Rn \subseteq \K(\Gamma)$.  Suppose that for each component $M_i$ of $\bdry W$ on which $\phi$ is nontrivial, $\rank_{\Z[\Gamma]}(H_1(M_i;\Z[\Gamma]))=\beta_1(M_i)-1$.  Then if $P$ is the kernel of the inclusion induced map $T(H_1(\bdry W;\Rn))\to T(H_1(W;\Rn))$ then $P$ is isotropic  with respect to the Blanchfield form on $T(H_1(\bdry W;\Rn)).$
\end{theorem*}

Now, for any PTFA group, $\Gamma$, $\Gamma^{(k)}=1$ for some $k$.  Since the cobordisms in which  we are interested are $(k)$-null bordisms for every $k$, this condition imposes only the restriction that $\Gamma$ be PTFA in our setting.  Since the components of $V_\infect$ and $V_\#$ are all given by zero surgery along knots The condition that $\rank_{\Z[\Gamma]}(H_1(M_i;\Z[\Gamma]))=\beta_1(M_i)-1=0$ holds for all coefficient systems which are nontrivial on $M_i$ by \cite[Proposition 3.10]{C}.  

In the case of $V_\#$, the meridian of any one of the components normally generates $\pi_1(V_\#)$ so that if $\phi$ is nontrivial on $V_\#$ then it is nontrivial on every boundary component.  In the case of $V_\infect$, the meridians of $R$ and $R_\eta(K)$ each normally generate $\pi_1(V_\infect)$ so that $\phi$ is  nontrivial on the $M(R)$ and $M(R_\eta(K))$ boundary components

\begin{proposition}\label{cobordism lemma}
Consider knots $K$ and $J$.  Let $\eta$ be an unknotted curve in the complement of $K$.  Let $V$ be either $V_\#(K,J)$ or $V_\infect(K,\eta,J)$ Given a nontrivial PTFA coefficient system $\phi:\pi_1(V_\#)\to\Gamma$ and $S\subseteq \Q[\Gamma]$ a right divisor set which is closed under the involution on $\Q[\Gamma]$.  Let $\Rn= \Q[\Gamma]S^{-1}$.  Then
\begin{enumerate}[label=(\arabic*)]
\item\label{cl sign} $\sigma^{(2)}(V,\Gamma)=\sigma(V)=0$,
\item\label{cl consum} $\ker\left(H_1\left(\bdry V_\#;\Rn\right) \to H_1\left(V_\#\Rn\right)\right)$ is isotropic with respect to $Bl_\Gamma$.
\item\label{cl infect1} If $V=V_\infect$ and $\phi$ is nontrivial on $M(J)$, then $$\ker\left(H_1\left(\bdry V_\infect;\Rn\right) \to H_1\left(V_\#;\Rn\right)\right)$$ is isotropic with respect to $Bl_\Gamma$.
\item\label{cl infect 2} If $\phi$ is trivial on $M(J)$, then
$$\ker\left(H_1\left(M(K(\eta, J));\Rn\right)\oplus H_1\left(M(K);\Rn\right) \to H_1\left(V_\infect;\Rn\right)\right)$$
 is isotropic with respect to $Bl_\Gamma$.
\end{enumerate}
\end{proposition}
\begin{proof}
  By \cite[Theorem 5.9]{blanchfieldDuality}, $\sigma^{(2)}(W,\Gamma)=\sigma(W)=0$.

The remaining conclusions are immediate consequences of \cite[Theorem 6.3]{blanchfieldDuality}.  Conclusions (2) and (3) follow since $H_2(\bdry V;\Rn)$ is torsion in these cases.  Conclusion (4) holds since if $\phi$ is trivial on $M(J)$, then the torsion part of $H_2(\bdry V_\infect;\Rn)$ is $H_1\left(M(K(\eta, J));\Rn\right)\oplus H_1\left(M(K);\Rn\right)$.  
\end{proof}

Finally, we prove the technical lemmas needed in the proof of Lemma~\ref{solvable to C}.  For convenience we restate them as we prove them.

\begin{replemma}{connected sum inheritance}
Let $K_i=K_{i,1} \# \dots \#K_{i,m_i}$ for $i=1, \dots, p$.  If $\iterate{\sqcup}{i=1}{p}M(K_i)$ bounds a $C(n,h)$, $(W,\Gamma)$, then $\iterate{\sqcup}{i=1}{p}\iterate{\sqcup}{j=1}{m_i}M(K_{i,j})$ bounds a $C(n,h)$.
\end{replemma}
\begin{proof}
Construct a new 4-manifold, $\widehat{W}$, by gluing to the 
$$M\left(K_i\right) = M\left(K_{i,1} \# \dots \#K_{i,m_i}\right)$$ 
boundary component of $W$ a copy of $V_\#\left(K_{i,1}, \dots, K_{i,m_i}\right)$ which we call $V_i$.  Do this for each $i$.  The resulting 4-manifold has boundary given by $\iterate{\sqcup}{i=1}{p}\iterate{\sqcup}{j=1}{m_i} M\left(K_{i,j}\right)$.  Since the map $i_*:H_1(M(K_i))\to H_1(V_i)$ is an isomorphism, one can use Condition \ref{C2} to conclude that $\phi:\pi_1(W)\to\Gamma$ extends over $\pi_1(V_i)$ for each $i$.  Specifically, if $\alpha_i$ is the monomorphism which exists since $(W,\Gamma)$ satisfies \ref{C2}, one can define the extension of $\phi$ to $\pi_1(V_i)$ by the composition:
\begin{equation*}
\pi_1(V_i)\to H_1(V_i) \overset{i_*^{-1}}{\to} H_1(M(K_i))\overset{\alpha_i}{\hookrightarrow} \Gamma
\end{equation*}

We claim that $(\widehat{W},\Gamma)$ is a $C(k,n)$.  Since the underlying group, $\Gamma$, did not change, Condition \ref{C1} still holds.  In order to see Condition \ref{C2} consider the following diagram:
\begin{equation*}
\begin{diagram}
\node{\pi_1(M(K_{i,j}))}\arrow{e}\arrow{s}\node{H_1(M(K_{i,j}))}\arrow{s,r}{\cong}\\
\node{\pi_1(V_\#)}\arrow{e}\node{H_1(V_\#)}\arrow{s,lr}{i_*^{-1}}{\cong}\\
\node{\pi_1(M(K_i))}\arrow{n}\arrow{e}\arrow{s}\node{H_1(M(K_{i}))}\arrow{s,r,J}{\alpha_i}\\
\node{\pi_1(W)} \arrow{e} \node{\Gamma}
\end{diagram}
\end{equation*}
The composition of the maps on the right hand column is the required monomorphism.  Its image is contained in the Abelian subgroup, $A$, given by the fact that ($W,\Gamma$) satisfies \ref{C2}.

We now use Proposition~\ref{cobordism lemma} \ref{cl consum}  to show Condition \ref{C3}.    Let $\psi:\pi_1(\widehat{W})\to\Lambda$ and $S$ be as in the statement of Condition \ref{C3}.  If $x$ and $y$ are in 
\begin{equation*}P:=\ker\left(\underset{i,j}{\oplus}H_1(M(K_{i,j});\Q[\Lambda]S^{-1})\to H_1(\widehat W;\Q[\Lambda]S^{-1})\right)\end{equation*}
  then there exist $x'$ and $y'$ in 
  \begin{equation*}P' := \ker\left(\underset{i}{\oplus}H_1(M(K_i);\Q[\Lambda]S^{-1})\to H_1(W;\Q[\Lambda]S^{-1})\right)\end{equation*}
   such that $x-x'$ and $y-y'$ are in $Q=\underset{i}\oplus Q_i$ where
   \begin{equation*}Q_i:=\ker\left( H_1(\bdry V_{i};\Q[\Lambda]S^{-1})\to H_1(V_{i};\Q[\Lambda]S^{-1})\right).\end{equation*}
Consider the equality from the sesquilinearity of the Blanchfield form, \begin{equation}\label{sesquilinear}Bl_\Lambda(x-x',y-y')=Bl_\Lambda(x,y)-Bl_\Lambda(x,y')-Bl_\Lambda(x',y)+Bl_\Lambda(x',y').\end{equation}
Since $Q$ is isotropic by Proposition~\ref{cobordism lemma} \ref{cl consum}, $Bl_\Lambda(x-x',y-y')=0$.  By assumption, Condition \ref{C3} holds for $(W,\Gamma)$, so $P'$ is isotropic and $Bl_\Lambda(x',y')=0$.  Since $x$ and $y'$ are carried by different components of $\bdry V$, as are $x'$ and $y$, $Bl_\Lambda(x,y')=Bl_\Lambda(y',x)=0$.  Thus, \eref{sesquilinear} reduces to $0=Bl_\Lambda(x,y)$ so that $P$ is isotropic and Condition \ref{C3} holds.

Condition \ref{C4} holds because of Novikov additivity, since by Proposition~\ref{cobordism lemma} \ref{cl sign}, $\sigma(V_i)=\sigma^{(2)}(V_i,\Theta)=0$.  This completes the proof.

\end{proof}

\begin{replemma}{infection inheritance}
  For $i=1, \dots m$, let $K_i$ be a knot, $R_i$ be slice and $\eta_i$ represent a doubly anisotropic element of $A_0(R_i)$.  If $\iterate{\sqcup}{i=1}{m}M(R_i(\eta_i,K_i))$ bounds a $C(n,h)$ with $n\le h$, $(W,\Gamma)$, then $\iterate{\sqcup}{i=1}{m}M(K_i)$ bounds a $C(n+1,h)$.
\end{replemma}
\begin{proof}

To each $M\left(R_i(\eta_i, K_i)\right)$-boundary component of $W$ glue a copy of $V_\infect (R_i,K_i,\eta_i)$ which we abbreviate as $V_i$.  Do this for each $i$. Call the resulting $4$-manifold $W_0$.  It has boundary $\bdry W_0=\iterate{\sqcup}{i}{} (M(K_i)\sqcup M(R_i))$.  For each $i$, let $E_i$ be the complement of a slice disk for the slice knot $R_i$.  Let $\widehat W$ be given by gluing $E_i$ to $W_0$ along the $M(R_i)$ boundary component for each $i$. 

Similar to the proof of Lemma~\ref{connected sum inheritance}, the coefficient system extends over $V_i\cup E$ and on $V_i\cup E$ it factors through Abelianization.  Since $\mu_i$ (the meridian of $K_i$) is isotopic in $V_i$ to $\eta_i$ which is nullhomologous, $\phi(\mu_i)$ is trivial.  Thus, $\mu_i\cong \eta_i$ lifts to a curve in the $\Gamma$-cover of $\widehat W$ and can be regarded as an element of $H_1(\widehat W; \Q[\Gamma])$.

If $x$ and $y$ are elements of 
$$P:=\ker\left(\underset{i}{\oplus}H_1(M(R_i);\Q[\Gamma])\to H_1(W_0;\Q[\Gamma])\right)$$
then there must exists $x'$ and $y'$ in 
$$P':=\ker\left(\underset{i}{\oplus}H_1(M(R_i(\eta_i,K_i));\Q[\Gamma])\to H_1(W;\Q[\Gamma])\right)$$
such that $x-x'$ and $y-y'$ are in $S=\underset{i}{\oplus}S_i$ where
$$S_i:=\ker\left(H_1(M(R_i);\Q[\Gamma])\oplus H_1(M(R_i(\eta_i,K_i));\Q[\Gamma])\to H_1(V_i;\Q[\Gamma])\right).$$
By Proposition~\ref{cobordism lemma} \ref{cl infect 2} $S$ is isotropic so that 
$$0=Bl(x-x',y-y')=Bl(x,y)-Bl(x,y')-Bl(x',y)+Bl(x',y').$$
By remark~\ref{C3 to isotropy}, $P'$ is isotropic so that $Bl(x',y')=0$.  Since $x$ and $y'$ as well as $x'$ and $y$ sit in different components, $Bl(x,y')=Bl(x',y)=0$.  Thus, $Bl(x,y)=0$ and $P$ is isotropic.

   Since the inclusion induced map $\nquotient{M(R_i)}{1}\to \nquotient{M(E_i)}{1}\cong \Z$ is an isomorphism, it follows that if $Q_i=\ker(A_0(R_i)\to A_0(E_i))$, then 
   $$
   \ker\left(H_1(M(R_i);\Q[\Gamma])
   \to H_1(E_i;\Q[\Gamma])
   \right)=Q_i\otimes \Q[\Gamma].
   $$
   
   By a Mayer-Vietoris argument, \begin{equation*}\ker\left(\underset{i}{\oplus}H_1(M(R_i);\Q[\Gamma])\to H_1(W;\Q[\Gamma])\right)=P+\underset{i}{\oplus}(Q_i\otimes \Q[\Gamma]).\end{equation*} 

If $\langle0,\dots,0,\eta_j\otimes1,0,\dots,0\rangle\in \underset{i}{\oplus}H_1(M(R_i);\Q[\Gamma])$ were in $P+\underset{i}{\oplus}(Q_i\otimes\Q[\Gamma])$, then there would exist some $p=\langle p_i \rangle\in P$ and $\langle q_i\rangle\in \underset{i}{\oplus} (Q_i\otimes\Q[\Gamma])$ with $\eta_j\otimes1 = p_j+q_j$ and $0=p_i+q_i$ when $i\neq j$.  Since $q_i$ is in the isotropic submodule $Q_i\otimes \Q[\Gamma]$ for each $i$ and $P$ is isotropic, this implies
\begin{equation*}
\begin{array}{rcl}
0&=&Bl_\Gamma(p,p)=\displaystyle \sum_i(Bl_\Gamma(p_i,p_i)) = Bl_\Gamma(p_j,p_j)+\sum_{i\neq j} Bl_\Gamma(q_i,q_i)\\ &=& Bl_\Gamma(p_j,p_j)
\end{array}\end{equation*}
So that $\eta_j\otimes 1=p_j+q_j$ sits in the sum of the isotropic submodule $Q_j\otimes \Q[\Gamma]$ together with the isotropic submodule generated by $p_j$.   Corollary~\ref{2xAnisoCor} then contradicts the assumption that $\eta_j$ be doubly anisotorpic.  Thus, $\mu_j\cong \eta_j$ must be nonzero in $H_1(\widehat W;\Q[\Gamma])$ and $H_1(M(K_i))$ maps injectively to $H_1(\widehat W;\Q[\Gamma])$.

Letting $G=\ker(\pi_1(\widehat W)\to \Gamma)$, define $\widecap\Gamma$ to be the quotient $\displaystyle\dfrac{\pi_1(W)}{G^{(1)}_\mathfrak{r}}$ where $G^{(1)}_\mathfrak{r}$ is the first term in the rational derived series of $G$.  Consider the following short exact sequence,
\begin{equation*}\displaystyle 0\to\dfrac{G}{G^{(1)}_\mathfrak{r}}\to\widecap{\Gamma}\to\Gamma\to0.\end{equation*}
It reveals first that $\widecap{\Gamma}$ is PTFA, since $\Gamma$ is PTFA and $\displaystyle \dfrac{G}{G^{(1)}_\mathfrak{r}}$ is TFA.  Secondly, since $G$ is the fundamental group of the $\Gamma$ cover of $\widecap W$, $\displaystyle\dfrac{G}{G^{(1)}_\mathfrak{r}}$ is the quotient of $H_1(\widecap W;\Z[\Gamma])$ by its $\Z$-torsion, into which $H_1(M(K_i))$ was shown to inject.  Thus, this choice of $(\widecap W, \widecap{\Gamma})$ satisfies \ref{C2}.  It satisfies \ref{C1} for $n+1$ since $\widehat{\Gamma}^{(n+1)}_\mathfrak{r}$ sits in the TFA group $\displaystyle \dfrac{G}{G^{(1)}_\mathfrak{r}}$, so that $\widecap\Gamma^{(k+2)}_\mathfrak{r}=0$.

The argument that $(\widehat{W}, \widehat{\Gamma})$ satisfies conditions \ref{C3} and \ref{C4} is just as in the proof of \ref{connected sum inheritance} with part \ref{cl infect1} of \ref{cobordism lemma} replacing part \ref{cl consum}.  
\end{proof}

\bibliographystyle{plain}
\bibliography{biblio}  

\end{document}